\documentclass[12pt]{amsart}
\usepackage{amssymb, amsmath}
\usepackage[mathscr]{euscript}
\usepackage{amscd}
\usepackage{stmaryrd}
\usepackage{pb-diagram}
\newcommand{\bbc}{\mathbb{C}}
\newcommand{\bbp}{\mathbb{P}}
\newcommand{\bbq}{\mathbb{Q}}
\newcommand{\bbr}{\mathbb{R}}
\newcommand{\bbz}{\mathbb{Z}}
\newcommand{\bbn}{\mathbb{N}}
\newcommand{\Sym}{\operatorname{Sym}}

\newtheorem{thm}{Theorem}[section]
\newtheorem{cor}[thm]{Corollary}
\newtheorem{lem}[thm]{Lemma}
\newtheorem{prop}[thm]{Proposition}

\newtheorem{defnn}[thm]{Definition}

\newtheorem{remarkk}[thm]{Remark}

\newtheorem{examplee}[thm]{Example}

\title[Mordell-Weil groups and the rank of elliptic curves over large fields]{ Mordell-Weil groups and the rank of
elliptic curves over large fields }
\author{Bo-Hae Im}
\date{January 20, 2003}
\address{Department of Mathematics, Indiana University, Bloomington,
Indiana 47405} \email{boim@indiana.edu} 

\begin{document}
\maketitle
\begin{abstract} Let $K$ be a  number field, $\overline{K}$ an algebraic closure of $K$ and $E/K$ an elliptic curve
defined over $K$. In this paper, we prove that if $E/K$ has a
$K$-rational point $P$ such that $2P\neq O$ and $3P\neq O$, then
for each $\sigma\in Gal(\overline{K}/K)$, the Mordell-Weil group
$E(\overline{K}^{\sigma})$ of $E$ over the fixed subfield of
$\overline{K}$ under $\sigma$ has infinite rank.
\end{abstract}

\section{Introduction}

 In \cite{fj74}, G.~Frey and M.~Jarden showed that if $K$ is an infinite field of finite type and $A$ is an abelian
variety of dimension $d\geq 1$ defined over $K$, then for any
positive integer $n$, there is a subset of $Gal(\overline{K}/K)^n$
of Haar measure 1 such that for every $n$-tuple
$(\sigma_1,\ldots,\sigma_n)$ belonging to the subset, the group of
rational points $A(\overline{K}(\sigma_1,\ldots,\sigma_n))$ of $A$
over the fixed subfield of $\overline{K}$ under
$(\sigma_1,\ldots,\sigma_n)$ has infinite rank.

In \cite{larsen}, M.~Larsen proved that for a number field $K$ and
an elliptic curve $E/K$ over $K$, there is a nonempty open subset
$\Sigma$ of $Gal(\overline{K}/K)$ such that for any $\sigma\in
\Sigma$, the Mordell-Weil group $E(\overline{K}^{\sigma})$ of $E$
over the fixed field under $\sigma$ has infinite rank.

It is natural to ask if such an open subset can be the whole
Galois group $Gal(\overline{K}/K)$. We have a positive answer for
elliptic curves defined over $\bbq$. In \cite{im}, we proved that
for any elliptic curve $E/\bbq$, the rank of
$E(\overline{\bbq}^{\sigma})$ is infinite, for every $\sigma\in
Gal(\overline{\bbq}/\bbq)$. Our approach in \cite{im} is
arithmetic: taking advantage of the modularity of elliptic curves
over $\bbq$ and the complex multiplication theory and constructing
an infinite supply of rational points of $E$ consisting of Heegner
points.

This paper is motivated by \cite{fj74}, \cite{im} and
\cite{larsen} and
 we prove in section 3 that if $E/K$ has a  $K$-rational point $P$ such that $2P\neq O$ and
 $3P\neq O$, then for each $\sigma \in
Gal(\overline{K}/K)$, the Mordell-Weil group
$E(\overline{K}^{\sigma})$ over the fixed subfield of
$\overline{K}$ under $\sigma$ has infinite rank. Here, we approach
by using Diophantine geometry which is a completely different
method from the one that we use in \cite{im}.

The main strategy for constructing infinitely many linearly
independent rational points of $E$ over $\overline{K}^{\sigma}$
for $\sigma\in Gal(\overline{K}/K)$ is approximately as follows:
 find a finite group $G$, a $\bbz$-free $\bbz[G]$-module $M$
and an infinite sequence $\{K_i/K\}_{i=1}^{\infty}$ of linearly
disjoint finite Galois extensions of $K$ with $Gal(K_i/K)\cong G$
such that for each $i$, $E(K_i)\otimes\bbq$  contains
 a $G$-submodule isomorphic to $M\otimes\bbq$. If $M^G=0$ but $M^g\neq 0$ for each $g\in G$, then we can find
 $\bbq$-independent points of $E(K_i)\cap E(\overline{K}^{\sigma})$ for any $\sigma\in Gal(\overline{K}/K)$.
For any pair $(G,M)$, $G$ acts on $E\otimes M\cong E^r$. Suppose
we find a projective line $\bbp^1$ in $(E\otimes M)/G$ over $K$.
If
 its preimage $X$  in $E\otimes M$ under the quotient map is an irreducible curve over $K$, then
 by the Hilbert irreducibility theorem (\cite{l83}, Chapter 9), most points in $\bbp^1(K)$ determine points in
$E^r(K_i)$ with $Gal(K_i/K)=G$; the coordinates generate the
desired $G$-submodule of $E(K_i)\otimes\bbq$. In this paper, we
take for $G$ the alternating group $A_n$ on $n=2k$ letters and for
the module $M$ the irreducible $(n-1)$-dimensional quotient of the
permutation representation of $A_n$.

In section 2, first, we show that $S_n$ admits a nontrivial action
on the $(n-1)$-fold product $E^{n-1}$ of $E$ and its quotient
$E^{n-1}/S_n$ by $S_n$ is isomorphic to the $(n-1)$-dimensional
projective space $\bbp^{n-1}$. We also find some properties of
transitive subgroups of $S_n$ which contain a transposition and
observe properties of subgroups of $A_n$ which occur as branched
Galois coverings of a projective line.

In section 3, if $K$ is totally imaginary and $E/K$ has a
$K$-rational point $P$ such that $2P\neq O$ and $3P\neq O$, then
we show that for some even integer $n$, there is a projective line
over $K$ in $E^{n-1}/S_n$ whose preimage in $E^{n-1}/A_n$ under
the double
 cover is a curve of genus 0, which gives infinitely many linearly independent points of $E$ over the
 fixed field of each
$\sigma\in Gal(\overline{K}/K)$.
%In section 4, by constructing an infinite sequence of linearly
%disjoint totally imaginary extensions  of degree $\leq 4$ over
%$K$, we show that if $K$ is a number field and $E/K$ is an
%elliptic curve over $K$, then the set $\{\sigma\in
%Gal(\overline{K}/K)\mid $ the rank of $E(\overline{K}^{\sigma})$
%is infinite$\}$ contains an open subset of Haar measure 1.

In section 4, as a special case, by the Hilbert irreducibility
theorem (\cite{l83}, Chapter 9) and the density of the Hilbert
sets over $\bbq$ in $\bbr$, we prove that if $K$ is a number field
and $K_{ab}$ is the maximal abelian extension of $K$, then for any
complex conjugation automorphism $\sigma\in Gal(\overline{K}/K)$,
the rank of $E((K_{ab})^{\sigma})$ is infinite. Hence, the rank of
$E(\overline{K}^{\sigma})$ is infinite.

Then, in section 5, we show that if  $\sigma\in
Gal(\overline{K}/K)$ is not a complex conjugation automorphism,
then, there is a totally imaginary finite extension of $K$ which
is fixed under $\sigma$. So by applying this to extend the ground
field to a totally imaginary extension for such automorphisms in
$Gal(\overline{K}/K)$ and combining the result of infinite rank of
the case of totally imaginary number fields, and the case of
complex conjugation automorphisms, we get a more general result
that if $K$ is an arbitrary number field and $E/K$ has a
$K$-rational point $P$ such that $2P\neq O$ and $3P\neq O$, then
for each $\sigma\in Gal(\overline{K}/K)$ the rank of
$E(\overline{K}^{\sigma})$ is infinite.

\vspace{0.2 in}
\begin{center}
\sc{Acknowledgements}
\end{center}
\vspace{0.1 in}

I wish to thank my thesis advisor, Michael Larsen for suggesting
 this problem, his guidance, valuable discussions and helpful
comments on this paper.

\section{Action of $S_n$ on $E^{n-1}$ and branched Galois coverings of $\bbp^1$ }

Let $n\geq 2$ be an integer. First, let $S_n$ be the symmetric
group on $n$ letters and $A_n$ the alternating subgroup of $S_n$.
Denote the $n$-fold product of $E$ by $E^n$. Naturally, $S_n$ acts
on $E^n$ by permutation, \emph{i.e.} if we denote its action by
`$\cdot$', for $\sigma \in S_n$ and an $n$-tuple
$(P_1,\ldots,P_n)\in E^n$,
$\sigma\cdot(P_1,\ldots,P_n)=(P_{\sigma(1)},\ldots,P_{\sigma(n)})$.
So does $A_n$ on $E^n$. Let $\Sigma :E^n \rightarrow E$ be the map
defined by the sum of coordinates of an $n$-tuple. Then identify
$E^{n-1}$ with $n$-tuples of elements in $E$ which sum to $O$
\emph{i.e.} $Ker(\Sigma)$. $S_n$ still acts on $E^{n-1}\cong
Ker(\Sigma)$ by the nontrivial induced permutation action.

Through the paper, we always consider $E^{n-1}$ as $Ker(\Sigma)$
so that a point in $E^{n-1}$ (or its quotient $E^{n-1}/S_n$ by
$S_n$) is a $n$-tuple $(P_1,\ldots,P_n)\in E^{n-1}$ whose
coordinates sum to $O$.

The following lemma gives the structure of the quotient space
$E^{n-1}/S_n$ of $E^{n-1}$ by $S_n$.

\begin{lem}\label{lem:action}
For each $n\geq 2$, $S_n$ admits a nontrivial action on $E^{n-1}$.
And the quotient space $E^{n-1}/S_n$ of $E^{n-1}$ by $S_n$ is
isomorphic to  the $(n-1)$-dimensional projective space
$\bbp^{n-1}$.
\end{lem}

\begin{proof}
Identify $E^{n-1}$ with $n$-tuples of elements in $E$ which sum to
$O$ \emph{i.e.} with the set $Ker(\Sigma)$, where $\Sigma :E^n
\rightarrow E$ is the map defined by the sum of coordinates of an
$n$-tuple. Then,  for each $(P_1,\ldots,P_n)\in Ker(\Sigma)\cong
E^{n-1}$, there is a rational function $f$ on $E$ such that
$\sum\limits _{i=1}^n (P_i)=(f)+n(O)$ as divisors. This gives a
map from $E^{n-1}$ to the linear space of all rational functions
$f$ on $E$ such that $(f)+n(O) \geq 0$. Denote this linear space
by $|n(O)|$.
%effective divisors on $E$ of degree $n$
%and linearly equivalent to the divisor $n(O)$ (which is called the
%complete linear system of $n(O)$ denoted by $|n(O)|$) given by
%$(P_1,\ldots,P_n)\mapsto (P_1)+\dots+(P_n)$. This is surjective and
%we can see that two $n$-tuples map onto the same divisor if and
%only if two $n$-tuples are the same up to permutations of $S_n$.
%This implies that the quotient space $E^{n-1}/S_n$ is isomorphic
%to $|n(O)|$.

Then, by the Riemann-Roch Theorem (\cite{h52}, Chapter IV, Theorem
1.3), the dimension of this space is $n$ as a vector space so it
gives an  ($n-1$)-dimensional projective space. We choose a basis
$f_0,\ldots, f_{n-1}$ of the space $|n(O)|$ and define a map $\phi
: E^{n-1} \rightarrow \bbp^{n-1}$ in the following way.

For each $(P_1,\ldots,P_n)\in Ker(\Sigma)\cong E^{n-1}$, there is
a rational function $f$ on $E$ such that $\sum\limits _{i=1}^n
(P_i)=(f)+n(O)$. Write $f =\sum\limits _{i=0}^{n-1}a_i f_i$ with
$a_0,\ldots,a_{n-1}\in\bbc$. Then, define
$\phi(P_1,\ldots,P_n)=(a_0:a_1:\cdots:a_{n-1})\in \bbp^{n-1}$.

Then two $n$-tuples which sum to $O$ in $E^{n-1}$ map onto the
same point in $\bbp^{n-1}$ under $\phi$ if and only if they are
the same up to permutations of $S_n$. This implies that the
quotient space $E^{n-1}/S_n$ is isomorphic to the projective space
$\bbp^{n-1}$.
 %and it
%is isomorphic to the $(n-1)$-dimensional projective space,
%$\bbp^{n-1}$. Therefore, the quotient space $E^{n-1}/S_n$ is
%isomorphic to $\bbp^{n-1}$.
\end{proof}

Now we find some properties of subgroups of $S_n$ which act
transitively on $\{1,2,\ldots,n\}$ and contain a transposition.
The following lemma assumes a weaker condition than in
(\cite{jar2}, Lemma 1.4).

\begin{lem}\label{lem:semidirect}
If $H$ is a subgroup of $S_n$ containing a transposition and $H$
acts transitively on $\{1,2,\ldots,n\}$, then there are positive
integers $m$ and $k$ such that  $mk=n$, where $m>1$, $k\geq 1$ and
there are a normal subgroup $M$ of $H$ and a subgroup $K$ of $S_k$
such that $M\cong (S_m)^k$, $H/M  \cong K$ and $K$ acts
transitively on $\{1,2,\ldots,k\}$, where $(S_m)^k=\underbrace{S_m
\times \cdot\cdot\cdot \times S_m}_{k~ times}$.

Moreover, $M\cap A_n \unlhd H\cap A_n$ and $(H\cap A_n)/(M\cap
A_n)\cong K$ and $H\cap A_n$ acts transitively on
$\{1,2,\ldots,n\}$. In particular, if $n$ is a prime $p$, then
$H\cong S_p$ and $H\cap A_p\cong A_p$.
\end{lem}

\begin{proof} Without loss of generality, we may assume that the transposition $(12)\in H$.
Define a relation $\sim$ on $\{1,2,\ldots,n\}$ by:  for $x,y\in
\{1,2,\ldots,n\}$,  $$x\sim y ~~\mbox{if and only if }~~ x=y
\mbox{ or there is a transposition }(xy)\in H.$$ Then, this
relation is an equivalence relation. In fact, the transitivity of
the relation holds, since if $(xy)\in H$ and $(yz)\in H$, then
$(xz)=(xy)(yz)(xy)\in H$.

 Since $(12)\in H$, $1\sim 2$, hence there is at least one nonempty equivalence class of $\{1,2,\ldots,n\}$.
 Suppose $x\sim y$. Then $(xy)\in H$. Let $x'\in \{1,2,\ldots,n\}$. Since $H$ acts transitively on $\{1,2,\ldots,n\}$,
there is $h\in H$ such that $h(x)=x'$. Now
$h(xy)h^{-1}=(h(x)~h(y))=(x'~h(y))$, which is in $H$. Hence
$h(x)=x'\sim h(y)$, \emph{i.e.} $x\sim y$ iff $h(x)\sim h(y)$.
Therefore, each equivalence class has the same number of elements.

Let $k$ be the number of equivalence classes and let $C_1,\ldots,
C_k$ be the equivalence classes of $\{1,2,\ldots,n\}$. And let
$m=\frac{n}{k}$. Each class $C_i$ has $m$ elements. Note that
$m\geq 2$, since $(12)\in H$.

For each $h\in H$, on each class $C_i$, $h(C_i)=C_{h_i}$, for some
$h_i\in\{1,2,\ldots,k\}$, since we have showed in the above that
$x\sim y$ iff $h(x)\sim h(y)$. And $h$ gives a bijection of $C_i$
and $C_{h_i}$. Hence we have a natural map $\phi_h :
\{C_i\}_{1\leq i\leq k}\rightarrow \{C_i\}_{1\leq i\leq k}$
defined by $\phi_h(C_i)=C_{h_i}$ where $i, h_i=1,2,\ldots,k$. This
shows that $\phi_h$ permutes equivalence classes $C_1,\ldots,C_k$.
Hence we get a permutation $\sigma_h\in S_k$ such that
$\sigma_h(i)=h_i$, where $i_k$ is given by $h(C_i)=C_{h_i}$.

So we can define a map $\pi : H\rightarrow S_k$ given by
$\pi(h)=\sigma_h$ defined as above. Then $\pi$ is a group
homomorphism, since $hh'(x)=h(h'(x))$, for $h,h'\in H$ and $x\in
\{1,2,\ldots,n\}$.

Let $K=Image(\pi)$ and $M=ker(\pi)$. Then $M\unlhd H$ and $K\leq
S_k$. Moreover, $K$ acts transitively on $\{1,2,\ldots,k\}$, since
$C_1\sqcup C_2\sqcup\cdot\cdot\cdot \sqcup C_k=\{1,2,\ldots,n\}$
and $H$ acts on $\{1,2,\ldots,n\}$ transitively.

Now we show that $M \cong \underbrace{S_m \times \cdot\cdot\cdot
\times S_m}_{k~ times}:=(S_m)^k$. Let $S(C_i)$ be the group of all
permutations on elements of $C_i$. For any $h\in M$, $h$ has a
decomposition, $h=h_1h_2\cdot\cdot\cdot h_k$, where each
permutation $h_i$ is a product of disjoint cycles in $S(C_i)$,
since $h$ is stable on each class. If $h,g\in M$, let
$h=h_1h_2\cdot\cdot\cdot h_k$ and $g=g_1g_2\cdot\cdot\cdot g_k$,
where $h_i, g_i\in S(C_i)$, then $hg=h_1g_1h_2g_2\cdot\cdot\cdot
h_kg_k$, since $h_i$ and $g_j$ are disjoint for $i\neq j$. Hence
we get an injective homomorphism $f:M\rightarrow S(C_1)\times
\cdot\cdot\cdot \times S(C_k)$ defined by $f(h)=(h_1,\ldots,h_k)$,
where $h=h_1h_2\cdot\cdot\cdot h_k$ and $h_i\in S(C_i)$. Since
$S(C_i)\cong S_m$ is generated by transpositions and for any
$x_i,y_i\in C_i$, there is a transposition $(x_iy_i)\in H$,
$f((x_1y_1)\cdot\cdot\cdot (x_ky_k))=((x_1y_1),\ldots ,(x_ky_k))$.
Hence $f$ is surjective. Therefore,  $M\cong S(C_1)\times
\cdot\cdot\cdot \times S(C_k)\cong (S_m)^k$. By the first
isomorphism theorem,  $H/M\cong K$.

From the above, we get a short exact sequence of groups,
$$1\longrightarrow (S_m)^k \stackrel{f}{\longrightarrow} H
\stackrel{\pi}{\longrightarrow} K \longrightarrow 1. $$ Now we
show that the following sequence is exact: $$1\longrightarrow
(S_m)^k\cap A_n \stackrel{f'}{\longrightarrow} H\cap A_n
\stackrel{\pi'}{\longrightarrow} K \longrightarrow 1, $$ where
$f'$ is the restriction of the inclusion $f$ to $(S_m)^k\cap A_n$.

First, it is obvious that $f'$ is injective, since $f$ is
injective. Moreover, since $ker(\pi)=Image(f)$, we have
$ker(\pi')=ker(\pi)\cap A_n=Image(f)\cap A_n=Image(f')$. So this
implies the exactness of the middle one. Now we need to show
\emph{$\pi'$} is surjective. For any $\sigma\in K$, there is $h\in
H$ such that $\pi(h)=\sigma$, since $\pi$ is surjective. If $h$ is
an even permutation, then $h\in H\cap A_n$ and
$\pi'(h)=\pi(h)=\sigma$. If $h$ is not even, then consider
$\sigma$ as a permutation of $\{C_1,\ldots,C_k\}$ as in the above.
Then there are two distinct integers $i$ and $j\in
\{1,2,\ldots,k\}$ such that $\sigma(C_i)=C_j$. Since $C_i$ has at
least two elements, there are two elements $a,b\in C_i$.
$\emph{i.e.}$ $a\sim b\in C_i$. Hence $(ab)\in H$. Moreover,
$(ab)\in Ker(\pi)$ from the construction. Hence $(ab)\circ h$ is
an even permutation and $\pi'((ab)\circ h)=\pi((ab)\circ
h)=\pi(h)=\sigma$. Hence $\pi'$ is surjective. Therefore, $M\cap
A_n\cong (S_m)^k\cap A_n=ker(f') \unlhd H\cap A_n$ and $(H\cap
A_n)/((S_m)^k\cap A_n)\cong K$.

Now we show that $H\cap A_n$ acts transitively on
$\{1,2,\ldots,n\}$. If $k=1$, then $m=n$, hence $H\cong S_n$.
Therefore, $H\cap A_n = A_n$, which acts transitively on
$\{1,2,\ldots,n\}$. Assume that $k\geq 2$. Let $a,b\in
\{1,2,\ldots,n\}$. We need to find an even permutation $\sigma \in
H$ such that $\sigma(a)=b$. If both  $a$ and $b$ are in the same
class $C_i$ for some $i\in \{1,2,\ldots,k\}$, then there is
$(ab)\in H$. Since $k\geq 2$, we choose two distinct elements $c$
and $d\in C_j$ for some $j\neq i$. Then if let $\sigma=(ab)(cd)$,
then $\sigma\in H\cap A_n$ and $\sigma(a)=b$.

If $a$ and $b$ are in distinct classes, say $a\in C_i$ and $b\in
C_j$ for $i\neq j$, then there is $\tau\in K$ such that
$\tau(i)=j$. Since $\tau$ is a bijection between $C_i$ and $C_j$,
there are $b'\in C_j$ and $a'\in C_i$ such that $\tau(a)=b'$ and
$\tau(a')=b$.  If $\tau$ is an even permutation, then let
$\sigma=(aa')(bb')\circ \tau$. Then $\sigma(a)=b$ and $\sigma\in
H\cap A_n$. If $\tau$ is odd, then let $\sigma=(bb')\circ\tau$.
Then $\sigma(a)=b$ and $\sigma\in H\cap A_n$. This completes the
proof.
\end{proof}

\begin{lem}\label{lem:invari}
If $H$ is a transitive subgroup of $S_n$  and $(V,\rho)$ is the
permutation representation of $S_n$, then the restriction of
$(V,\rho)$ to $H$ has  one 1-dimensional invariant subspace.
\end{lem}

\begin{proof} Let $e_1,\ldots,e_n$ be a basis for the restriction  $(V,\rho')$ of the permutation representation of
$S_n$ to $H$. Let $H_1=\{h\in H \mid h(1)=1\}$ be the stabilizer
of 1 in $H$. Let $W$ be the subspace of $V$ generated by $e_1$.
Then $W$ is invariant under $H_1$. Moreover, since $H$ acts
transitively on $\{1,2,\ldots,n\}$, we have exactly $n$ left
cosets of $H_1$ in $H$. Hence we can identify the permutation
representation $(V,\rho')$ with the induced representation
$\left(\bigoplus \limits _{i=1}^n \rho'_{g_i}(W),~
\mbox{Ind}_{H_1}^H(1)\right)$ of $H$ by the trivial representation
$(W, 1)$ of $H_1$, where $g_i$ are representatives of left cosets
of $H_1$ in $H$.

If we denote  by $1$   the trivial representation of $H$, then by
Frobenius reciprocity, $$\langle ~1,
\mbox{Ind}_{H_1}^{H}(1)~\rangle _{H}
=\langle~\mbox{Res}^{H_1}_{H}(1),1~\rangle _{H_1}=\langle ~1,1~
\rangle _{H_1}=1.$$ Therefore, the restriction $(V,\rho)$ of the
permutation representation to $H$ has  one 1-dimensional invariant
subspace.
\end{proof}

\begin{cor}\label{cor:HA_n}
If $H$ is a transitive subgroup of $S_n$ and $H$ contains a
transposition, then the restriction of the permutation
representation of $S_n$ to $H\cap A_n$ has one 1-dimensional
invariant subspace.
\end{cor}

\begin{proof} This follows from Lemma \ref{lem:invari}, since  $H\cap A_n$ acts transitively on $\{1,2,\ldots,n\}$ by
Lemma \ref{lem:semidirect}.
\end{proof}

\begin{lem}\label{lem:inva}
Let $n\geq 1$ be an integer. Let $\sigma \in A_n$ have $k$
disjoint cycles, for some positive integer $k\leq n$. Then there
are $k$ fixed vectors under $\sigma $ in the permutation
representation of $A_n$.
\end{lem}

\begin{proof} Let $e_1,\ldots, e_n$ be a basis of the permutation representation of $A_n$. Let $\sigma\in A_n$ have $k$
disjoint cycles. Then, they form $k$ partitions $C_1,\ldots,C_k$
of $\{1,2,\ldots,n\}$.

For $1\leq i\leq k$, let

$$v_i=\sum\limits_{j\in C_i} e_j.$$ Then, these
$k$ vectors are fixed under $\sigma$.
\end{proof}

\begin{lem}\label{lem:cycle} For any even integer $n$, every element in $A_n$ has more than one cycle.
\end{lem}

\begin{proof}  Let $\sigma\in A_n$ have the cycle decomposition
$\sigma = (a_{11}\cdot\cdot\cdot a_{1 m_1})(a_{21}\cdot\cdot\cdot
a_{2 m_2})\cdot\cdot\cdot(a_{k1}\cdot\cdot\cdot a_{k m_k})$, where
$a_{ij}\in\{1,2,\ldots,n\}$ are distinct and
$m_1+m_2+\cdot\cdot\cdot+m_k=n$ for some positive integer $m_i$.
If $k=1$, then $m_1=n$ and $\sigma=(a_{11}a_{12}\cdot\cdot\cdot
a_{1n})$ has one cycle of length of the even integer $n$ which is
an odd permutation, hence it is not in $A_n$. Thus, we must have
that $k\geq 2$. This implies that $\sigma\in A_n$ has at least two
cycles.
\end{proof}

The following two lemmas show that subgroups of $S_n$ which occur
as Galois coverings of a projective line (which is isomorphic to a
projective closure of a  base-point free linear system of $E$) act
transitively on $\{1,2,\ldots,n\}$ and contain a transposition.

\begin{lem}\label{lem:tran} Let $K$ be a number field and $E/K$ an elliptic curve over $K$.
Suppose that there is a projective line $L$ in $E^{n-1}/S_n \cong
\bbp^{n-1}$ which is a projective closure of a base point-free
linear system of $E$. Let $C$ be  the preimage of $L$ in
$E^{n-1}/A_n$  under the double cover from $E^{n-1}/A_n$ to
$E^{n-1}/S_n$ and let the preimage in $E^{n-1}$ of $C$ under the
quotient map of $E^{n-1}$ by $A_n$ have a decomposition $X_1\cup
X_2\cup\cdot\cdot\cdot \cup X_k$ into irreducible components
$X_i$.

Then, for each $i=1,\ldots,k$, the morphism $X_i\rightarrow C$ is
a Galois covering with $H_i=Gal(X_i/C)$ a subgroup of $A_n$ and
each $X_i\rightarrow L$ is also a Galois covering with
$M_i=Gal(X_i/L)$ such that $M_i\leq S_n$ and $M_i\cap A_n = H_i$.
Moreover, the $M_i$ are conjugate to each other and each $M_i$
acts transitively on $\{1,2,\ldots,n\}$.
\end{lem}

\begin{proof} For each $i=1,\ldots,k$, the morphism $\psi_i:X_i\rightarrow
C$ is the quotient map by the stabilizer $H_i$ of $X_i$  in $A_n$,
that is, $H_i=\{\sigma\in A_n \mid \sigma\cdot X_i=X_i\}$. Since
$X_i$ is an irreducible component of the preimage of $C$ under the
action of $A_n$, for any $\sigma\neq 1\in H_i$, $X_i$ is not
contained in the kernel of $1-\sigma$ acting on $E^n$. So $\psi_i$
is a regular Galois covering map with Galois group
$Gal(X_i/C)=H_i$ a subgroup of $A_n$. Similarly, each map
$\phi_i:X_i \rightarrow L$ is also a Galois covering with Galois
group $Gal(X_i/L)$ which is
 the stabilizer of $X_i$ in $S_n$ and $\phi_i$ is
the composite of the Galois covering from $X_i$ to $C$ with
$H_i=Gal(X_i/C)$ and the double cover from $C$ to $L$. Hence $M_i
\leq S_n$ and $M_i\cap A_n = H_i$.

First, we show that the $M_i$ are conjugate to each other. Note
that $S_n$ acts transitively on $\{X_1,X_2,\ldots,X_k\}$, since
$X_1\cup X_2\cup\cdot\cdot\cdot \cup X_k$ is the preimage of the
irreducible curve $L$. Hence for each $i$, there is $\tau_i\in
S_n$ such that $\tau_i\cdot X_i=X_1$. Let $\sigma\in M_1$. Then
$\tau_i\cdot X_i=X_1=\sigma\cdot X_1=\sigma\tau_i\cdot X_i$. Hence
$\tau_i^{-1}\sigma\tau_i\cdot X_i=X_i$. Hence
$\tau_i^{-1}\sigma\tau_i\in M_i$. This proves that
$\tau_i^{-1}M_1\tau_i=M_i$ for each $i$.

Next, we show that each $M_i$ acts transitively on
$\{1,2,\ldots,n\}$.  It is enough to show that $M_1$ acts
transitively on $\{1,2,\ldots,n\}$, since $M_i$ are conjugate to
each other. Note that $E^{n-1}/S_n$ is isomorphic to the
projective closure of the linear space $|n(O)|$ of all rational
functions $f$ such that $(f)+n(O)\geq 0$ by Lemma
\ref{lem:action}. Since the curve $L\subset E^{n-1}/S_n$ is a
projective closure of a base point-free linear system of $E$,
there exists an elliptic function $f$ in
$H^0(E,\mathcal{L}(n(O)))$ which has $n$ zeros which sum to
 $O$ of $E$ such that the base-point free linear system is
generated by $f$ and the constant function $1$.
%Note that
%the map $\phi_1:X_1\rightarrow L$ is surjective. In fact,
%if it is not surjective, then it is a constant map. Then
%since they are conjugate to each other from the above,
%every map $\phi_i$, for $i=1,2,\ldots,k$ is also a constant
%map. This is impossible, since the map $X_1\cup
%X_2\cup\cdot\cdot\cdot \cup X_k \rightarrow C$ is
%surjective.

 Parameterize an open dense subset of the projective line $L$ in $\bbp^{n-1}\cong E^{n-1}/S_n$
  by the parameter $\lambda$ such that
 $f-\lambda$ represents a point of the open subset. Let
 $g:E^{n-1} \rightarrow E$ be defined by
 $g(P_1,\ldots,P_n)=P_1$ where the sum of $(P_1,\ldots,P_n)$ is
 $O$ as a point of $E^{n-1}$. Then,
 the curve $X_1 \subseteq E^{n-1}$ maps to $E$ through $g$
 as well as to the projective line $L\in E^{n-1}/S_n$
 through the quotient map by $S_n$. So $X_1$ maps to
 $E\times L$ and projects onto $L$.

Choose a fundamental domain so that the distinct zeros
$z_1,\ldots,z_n$ of $f-\lambda$ are in the interior of the domain.
Let $i\in \{2,3,\ldots,n\}$ be fixed. We can take a path from
$z_1$ to $z_i$ so that the path doesn't pass through the other
zeros of $f-\lambda$. Moving along this path in $E^{n-1}$, let
$f-\lambda '$ be the image in $E^{n-1}/S_n$ of the end point of
the path. Then $z_i$ is a common zero of both $f-\lambda$ and
$f-\lambda '$, we get $\lambda=\lambda'$. Hence  the image of the
path in $E^{n-1}/S_n$ is a closed loop starting and ending at
$\lambda$, that is, $f-\lambda$. This implies that the path from
$z_1$ to $z_i$ stays in the connected component $X_1$. Since
$\phi_1$ is a Galois covering with the Galois group $M_1$ which is
the stabilizer of $X_1$ in $S_n$ and $X_1$ maps to $E\times L$
through $g$ which is defined as a first coordinate of a point of
the preimage of $f-\lambda\in L$ and the morphism $\phi_1$, and
the first coordinates of the preimage of $f-\lambda$ are $z_1$ and
$z_i$ under the closed loop, this shows that there is an element
$\sigma\in M_1$ such that $\sigma(1)=i$. This shows $M_1$ acts
transitively on $\{1,2,\ldots,n\}$.
\end{proof}

The following lemma has a similar setting as in (\cite{jar2},
Lemma 1.5). But here, we assume that there is a divisor in a given
projective line which decomposes into the sum of one ramified
divisor of degree 2 and other divisors of odd degree or unramified
divisors under a Galois covering, while the lemma in (\cite{jar2},
Lemma 1.5) assumes every divisor decomposes into one ramified
divisor of degree 2 and other unramified divisors.

\begin{lem}\label{lem:(12)}
Suppose there is a curve $L \subset E^{n-1}/S_n\cong \bbp^{n-1}$
which is isomorphic to a projective closure of a base-point free
linear system on $E$ and the normalization of its preimage in
$E^{n-1}$ under the quotient map is $X_1\cup
X_2\cup\cdot\cdot\cdot \cup X_k$ such that for each
$m\in\{1,2,\ldots,k\}$,
 the Galois covering $M_m:=Gal(X_m/L)$ is a subgroup of $S_n$.
 Then if $L$ contains a divisor $D=2(P_1)+\sum\limits _{i=2}^{\ell} k_i (P_i) -n(O)$, where  $P_i$ are points
 of $E$ such
 that $P_i\neq P_j$, $2P_1+\sum\limits _{i=2}^{\ell} k_i P_i=O$ and $k_i$ are odd integers $\geq 1$ with $\sum\limits
 _{i=2}^{\ell} k_i=n-2$,
then  each $M_m$ contains a transposition.
\end{lem}

\begin{proof} Note that
if  $k_i=1$, for all $ i=2,\ldots,\ell$, then we apply the proof
in (\cite{jar2}, Lemma 1.5) with the given divisor $D$ to get a
transposition. Now we assume the general case when $k_i$ are odd
integers.

 For each $m=1,\ldots,k$, let
$\phi_i: X_i\rightarrow L$ be the restriction of the quotient map
of $E^{n-1}$ by $S_n$ with $M_m=Gal(X_m/L)$. Let $M_m$ act by
permutation of coordinates of each point: for $\sigma\in M_m$,
$\sigma\cdot
(P_1,P_2,\ldots,P_n)=(P_{\sigma(1)},P_{\sigma(2)},\ldots,P_{\sigma(n)})$,
where $P_n=-(P_1+P_2+\cdot\cdot\cdot+P_{n-1})$.

Suppose $L$ contains a divisor $D=2(P_1)+\sum\limits _{i=2}^{\ell}
k_i (P_i) -n(O)$, where $P_1,\ldots,P_{\ell}$ are distinct points
of $E$ and $k_i$ are odd integers such that
$\sum\limits_{i=2}^{\ell} k_i=n-2$. Let $f$ be the function whose
divisor is equivalent to $D$ and let $z_i$ be the zeros of $f$
corresponding to $P_i$ for each $i=1,\ldots,{\ell}$, \emph{i.e.}
$f(z)=(z-z_i)^{k_i}\left(a_{i0}+a_{i1}(z-z_i)+\cdot\cdot\cdot+\right)$
with $a_{i0}\neq 0$. Then by Hensel's Lemma (\cite{sil86}, Chapter
IV, Lemma 1.2), for a number $\lambda$ with small $|\lambda|$,
 $f-\lambda=(z-z_i)^{k_i}(a_{i0}+a_{i1}(z-z_i)+\cdot\cdot\cdot+\cdot\cdot)-\lambda$ has zeros at
$$Q_{1,1}=z_1+\left(\frac{\lambda}{a_{10}}\right)^{\frac{1}{2}}+A_1,
~~~
Q_{1,2}=z_1-\left(\frac{\lambda}{a_{10}}\right)^{\frac{1}{2}}+A_2,$$
$$\mbox{and }
Q_{i,j}=z_i+\left(\frac{\lambda}{a_{i0}}\right)^{\frac{1}{k_i}}\zeta
_{k_i}^{j-1}+B_{i,j},  \mbox{ for each } i=2,\ldots,{\ell}, \mbox{
and } j=1,\ldots,k_i,$$ where $A_1,A_2,$ and $B_{i,j}$ are
convergent Puiseux series in $\lambda$ such that each term of
$A_1$ and $A_2$ is of higher degree than $\lambda^{\frac{1}{2}}$
and each term of $Q_{i,j}$ is of higher degree than
$\lambda^{\frac{1}{k_i}}$, and $\zeta_{k_i}$ is a primitive
$k_i$-th root of unity.

Note that each quotient map $\phi_m$ by $M_m$ is surjective. Let
$k_1=2$.  Choose a small enough number $\lambda$ such that for all
$i=1,\ldots,{\ell}$, the circles centered at $z_i$ with radius
$\left(\frac{\lambda}{a_{i0}}\right)^{\frac{1}{k_i}}$ do not
intersect each other. For each $i=1,\ldots,{\ell}$, $k_i$ points
$Q_{i,j}$ for $j=1,\ldots,k_i$ lie in the circle centered at
$z_i$.

Since these circles are closed paths and $Q_{i,j}$ are zeros of
one fixed function $f-\lambda$,  their preimages in $E^{n-1}$
still stay in one component $X_m$ for some $m$. Therefore, for
each $i=1,\ldots,{\ell}$, there exists a cycle $\tau_i$ in $S_n$
of length $k_i$ which permutes $Q_{i,j}$ for $j=1,\ldots,k_i$ and
the product of all $\tau_i$ is in $M_m$. In particular, $\tau_i$
are disjoint cycles and $\tau_1$ is a transposition permuting
$Q_{1,1}$ and $Q_{1,2}$.

Let $c=lcm(k_2,\ldots,k_{\ell})$. Then, $c$ is odd, since all
$k_i$ for $i=2,\ldots,{\ell}$ are odd. Since $\tau_1$ is of order
$2$, $(\tau_1\tau_2\cdot\cdot\cdot\tau_{\ell})^c=\tau_1\in M_m$.
Hence, $M_m$ contains a transposition $\tau_1$. Since $M_m$ are
conjugate to each other by Lemma \ref{lem:tran}, every $M_m$ has a
transposition.
\end{proof}

\section{$E$ over totally imaginary number fields with a rational point $P$ such that
$2P\neq O$ and $3P\neq O$}

First, we show that if $K$ is a totally imaginary number field and
$E/K$ has a $K$-rational point $P$ such that $6P\neq O$, then, for
some even integer $n$, there is a projective line over $K$  in
$E^{n-1}/A_n \cong \bbp^{n-1}$ whose preimage under the quotient
map of $E^{n-1}$ by $A_n$ is a curve of genus 0 in $E^{n-1}$ over
$K$. We will need the following lemma to show the existence of
such a projective line. We start with the definition of
\emph{rank} of quadratic forms that we use in this paper.

\begin{defnn} The   rank   of a quadratic form $\Phi$   on the space $V$ is
the codimension of the orthogonal complement of   $V$   with
respect to   $\Phi$   in the sense of $($\cite{ser}, Chapter IV,
Section 1.2, pp.28$)$.
\end{defnn}

\begin{lem}\label{lem:intersec} Suppose $\Phi_1$ and $\Phi_2$ are
two quadratic forms defined over $\overline{K}$  such that for all
$r,s\in \overline{K}$, not both zero, the form $r\Phi_1+s\Phi_2$
is of rank $\geq 5$. Then the intersection of the zero loci of
$\Phi_1$ and $\Phi_2$ is not entirely contained in a finite union
of hyperplanes.
\end{lem}
\begin{proof} By the abuse of the notation, we denote the intersection of two
hypersurfaces defined by $f$ and $g$ by $f\cap g$ and the union of
them by $f\cup g$.

Suppose codim$( \Phi_1 \cap\Phi_2\cap L)=2$ for some hyperplane by
$L$. If both intersections $\Phi_1\cap L$ and $\Phi_2\cap L$ are
irreducible, then, $\Phi_1\cap L=\Phi_2\cap L$. Thus, $\Phi_1
\equiv c\Phi_2 (\bmod~ L)$, for some $c\in \overline{K}$, that is,
$\Phi_1-c\Phi_2=LL'$ for some linear form $L'$. Hence, the pencil
of $\Phi_1$ and $\Phi_2$ contains some form which has rank $\leq
2$, which leads to a contradiction to the hypothesis.

So we assume that a quadratic form, say $\Phi_1$ intersected with
$L$ is reducible into two hyperplanes defined by linear forms
$L_2$ and $L_3$ on the original space. Then, $\Phi_1\equiv L_2L_3
~(\bmod~ L)$. Therefore, for some linear form $L_4$,
$\Phi_1=LL_4+L_2L_3$ so it has rank $\leq 4$, which is a
contradiction to the hypothesis. Hence, we have shown that for
every hyperplane by $L$,
$$\mbox{codim}( \Phi_1 \cap\Phi_2\cap L)< 2.$$

 Now, suppose the intersection of $\Phi_1$ and
$\Phi_2$ is entirely contained in the union of hyperplanes
 by $L_1, \ldots,L_n$. Then, $$\min\limits_{1\leq i\leq n}\mbox{codim}( \Phi_1\cap\Phi_2
\cap L_i)=\mbox{codim}(\Phi_1\cap\Phi_2)=2,$$ which is impossible.
This completes the proof.
\end{proof}

We will need the following weak approximation of quadrics.

\begin{prop}
\label{prop:appro} There exists a function
$F:\;\bbn\rightarrow\bbn$ with the following property: Given a
non-negative integer $n$, a number field $K$, a $K$-vector space
$V$, an $n$-dimensional $K$-vector space of quadratic forms
$W\subset \Sym^2 V$ on $V$, and a finite set of places $S$ of $K$,
if for every non-zero $w\in W$, there exists an $F(n)$-dimensional
subspace $V_w\subset V$ on which $w$ is non-degenerate, then the
intersection of all quadrics in $W$, $X_W(K)$, is dense in
$\prod_{v\in S} X_W(K_s)$.
\end{prop}
\begin{proof} See (\cite{imlarsen}, Theorem 5).
\end{proof}

\begin{prop}
\label{prop:An} Let $K$ be a totally imaginary number field. If
$E/K$ has a  $K$-rational point $P$ such that $2P\neq O$ and
$3P\neq O$,
 then for some even integer $n$,  there is a projective line over $K$ in
 $E^{n-1}/S_n \cong \bbp^{n-1}$ as a projective closure of
 a  base-point free linear system of $E$
 such that the normalization of its preimage under the double cover is a curve of
genus 0 over $K$ in $E^{n-1}/A_n$ which contains a divisor of a
rational function on $E$ of degree $n$ which has one double zero
and all other zeros of odd order (including simple zeros).
\end{prop}

\begin{proof} By Lemma \ref{lem:action}, $E^{n-1}/S_n\cong \bbp^{n-1}$ which is isomorphic
to the $(n-1)$-dimensional projective space
$\bbp(H^0(E,\mathcal{L}(n(O))))$.

If $f$ is an elliptic function of degree $n$, holomorphic except
at a unique pole $O$, the vector space spanned by $f$ and $1$
defines a pencil of all divisors $(a+bf)+n(O)$ with $a,b\in\bbc$
on $E$ linearly equivalent to $n(O)$, or equivalently, a line on
$E^{n-1}/S_n\cong \bbp^{n-1}$.
 Note that since $P$ is neither 2-torsion nor
3-torsion, we have that $-2P\notin \{P, O\}$.
 Now we find an elliptic
function $f$ of degree $n=2k$ for some integer $k$, whose
derivative is of the form $f'=lh^2$, where
 $l$ has the divisor
$2(P)+(-2P)-3(O)$ and $h$ is in the vector space of elliptic
functions defined over $K$ with divisors $\geq (1-k)(O)$.
 Let $y+ax+b=0$ be the affine tangent line at a $K$-rational
point $P$ and let $l:=y+ax+b$. Let $f$ and $f'$ be as follows:

\vspace{.5cm}

 Case I : Suppose $n\equiv 0$ (mod 4). Let
$n=4m$ for some integer $m$. For parameters $a_0,\ldots,a_{m-1}$,
$b_0,\ldots,b_{m-3}$, $d_0,\ldots,d_{2m-2}$, $c_1,\ldots,c_{2m}$
to be determined and the given tangent line $l=0$, let $$f(z)=
y(d_{2m-2}x^{2m-2}+\cdot\cdot\cdot+d_1x+d_0)+c_{2m}x^{2m}+\cdot\cdot\cdot+c_1x
$$$$\mbox{ and } f'(z)=l(h(z))^2, $$ where $ h(z)=
a_{m-1}x^{m-1}+\cdot\cdot\cdot+a_1x+a_0+y(x^{m-2}+\cdot\cdot\cdot+b_1x+b_0).$

\vspace{.5cm}
 Case II : Suppose $n\equiv 2$ (mod 4). Let
$n=4m+2$ for some integer $m$. For parameters
$a_0,\ldots,a_{m-1}$, $b_0,\ldots,b_{m-2}$, $d_0,\ldots,d_{2m-1}$,
$c_1,\ldots,c_{2m+1}$ to be determined and the given tangent line
$l=0$, let $$ f(z)=
y(d_{2m-1}x^{2m-1}+\cdot\cdot\cdot+d_1x+d_0)+c_{2m+1}x^{2m+1}+\cdot\cdot\cdot+c_1x
$$$$\mbox{ and } f'(z)=l(h(z))^2, $$ where
$h(z)=x^m+a_{m-1}x^{m-1}+\cdot\cdot\cdot+a_1x+a_0+y(b_{m-2}x^{m-2}+\cdot\cdot\cdot+b_
1x+b_0).$

\vspace{.5 cm}

 From the equations obtained by  equating the coefficient
of each $x^{i}y^j$-term of $f'(z)$ with that of the derivative of
$f(z)$ given in the above (equivalently, by equating $f$ with the
integral of $f'$ along two periods of $E$), we get two quadratic
equations over $K$ in
 $\frac{n-4}{2}$  variables, namely
$a_0,\ldots,a_{m-1}, b_0,\ldots, b_{m-4}$ and $b_{m-3}$ if $n=4m$,
and $a_0,\ldots,a_{m-1}, b_0,\ldots,b_{m-3}$ and $b_{m-2}$ if
$n=4m+2$. Homogenize these two quadratic equations to get two
quadratic forms in $\frac{n-4}{2}+1$ variables with a new
variable. We need to find a common isotropic vector over $K$ of
two quadratic forms which defines a common solutions of two
original quadratic equations, (that is, which is outside the
hyperplane at $\infty$) and defines $f'=lh^2$ such that
$h(-2P)\neq 0$ in the above notation of cases I and II.

 Let $D$ be a
fundamental domain of $E$ and $C_1$ and $C_2$ be two line segments
dividing the fundamental domain of $E$ into four congruent
parallelograms and  $I_1 $ and $I_2$  the first half line segments
of $C_1$ and $C_2$ respectively as shown below.

\vspace{.8cm}

\begin{center}
\setlength{\unitlength}{.8cm}
\begin{picture}(8,8)
  \put(-1,3){\circle*{.17}}
  \put(-1.7,2.3){  $O$}
  \put(-1,3){\line(1,0){6}}
  \put(0,5){\line(1,0){6}}
  \put(1,7){\line(1,0){6}}
  \put(-1,3){\line(1,2){2}}
  \put(2,3){\line(1,2){2}}
  \put(5,3){\line(1,2){2}}
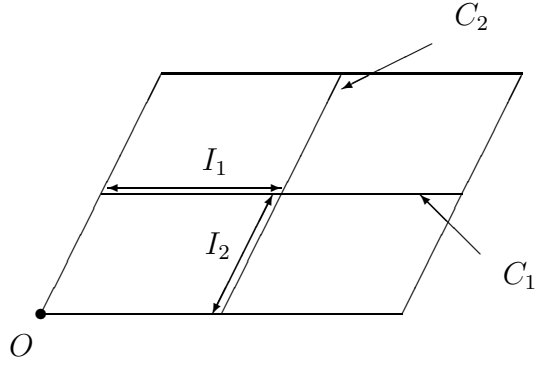
  \put(-3.5,.8){$\langle$A fundamental domain $D$ of $E$ with two periods
 $C_1$ and $C_2\rangle$}
   \put(1.5,5.4){ $I_1$}
   \put(0.3,5.1){\vector(1,0){2.7}}
    \put(2.8,5.1){\vector(-1,0){2.7}}
   \put(6.3,4){\vector(-1,1){1}}
   \put(6.5,3.5){ $C_1$}
    \put(1.55,4){ $I_2$}
    \put(1.9,3.1){\vector(1,2){.95}}
    \put(2.8,4.9){\vector(-1,-2){.95}}
    \put(5.5,7.5){\vector(-2,-1){1.5}}
    \put(5.7,7.8){  $C_2$}
\end{picture}
% \begin{pspicture}(4,4)
% \psset{linewidth=.5 cm}
%% \rput(-.27,-.27){$o$}
% \qdisk(0,0){2pt}
% \psline{-}(0,0)(4,0)
% \psline{-}(0,0)(1,3)
% \psline{-,linewidth=.05 cm}(2,0)(3,3)
% \psline{-,linewidth=.05 cm}(.5,1.5)(4.5,1.5)
% \psline[linestyle=dashed,dash=3pt 2pt](4,0)(5,3)
% \psline[linestyle=dashed,dash=3pt 2pt](1,3)(5,3)
% \psline{->}(4,3.6)(2.9,2.5)
% \rput(4.45,3.75){$c_2$}
% \psline{->}(5,.5)(4,1.47)
% \rput(5.2,.3){$c_1$}
% \psline{<->}(.55,1.6)(2.5,1.6)
% \psline{<->}(2.4,1.5)(1.9,0)
% \rput(1.5,1.9){$i_1$}
% \rput(1.7,.78){$i_2$}
% \rput(2.3,-1){\bf \tiny $\langle$a fundamental domain $d$ of $e$ with two periods
% $c_1$ and $c_2\rangle$}
%\end{pspicture}
\end{center}

\vspace{.5cm}

 Let $M=max\{F(2), 5\}$, where $F$ is the
function given in Proposition \ref{prop:appro}. We can choose $2M$
holomorphic functions $f_1,\ldots,f_{2M}$ on $I_1\cup I_2$ such
that
$$ \displaystyle\int_{I_1} lf_i f_j dz = \displaystyle\int_{I_2}
lf_i f_j dz =0, \mbox{ for  } i\neq j, \hspace{1.2 cm}$$
$$\hspace{.23 cm} \displaystyle\int_{I_1\cup I_2} lf_i ^2 dz =
\int _{I_1} lf_i^2 dz \neq 0, \mbox{ for } i=1,2,\ldots,M,$$ $$
\mbox{ and } \displaystyle\int_{I_1\cup I_2} lf_i ^2 dz = \int
_{I_2} lf_i^2 dz\neq 0, \mbox{ for }i=M+1,,\ldots,2M.$$

Since the Weierstrass $\wp$-function $x=\wp(z): I_1\cup
I_2\rightarrow \bbc$ is injective,  its inverse $\wp^{-1}$ is
well-defined on the image $\wp(I_1\cup I_2)$ and the image is a
compact contractible set in $\bbc$. Hence the complement of
$\wp(I_1\cup I_2)$ is connected.  So by Mergelyan's Theorem
(\cite{rudin}, pp.390), each holomorphic function $f_i\circ
\wp^{-1} :\wp(I_1\cup I_2)\rightarrow \bbc$ can be approximated by
some polynomial $p_i(z)$, for each $i=1,\dots, 2M$. Moreover,
since $K$ is a totally imaginary number field, $K$ is dense in
$\left(\prod\limits_{v\in S_{\infty}} \bbc\right)$  with respect
to the usual topology for any embeddings of $K$ in $\bbc$, where
$S_{\infty}$ is the set of all infinite places. Hence we may
assume that coefficients of $p_i(z)$ are in $K$. So each $f_i$ can
be approximated by the polynomial $p_i(x)$ in terms of $x=\wp(z)$
with coefficients in $K$.

Let $W$ be a space of dimensional $\geq 2M$ generated by elliptic
functions including all of $p_i(x)$ for $i=1,\ldots,2M$. Then, any
two quadratic forms $\Phi_1$ and $\Phi_2$ over $K$ obtained from
the homogenization of the integration on $W$ satisfy the property:
$$\mbox{for any } r,s\in \overline{K} \mbox{ not both zero, any form in the pencil }
r\Phi_1+s\Phi_2 \mbox{ is of rank }\geq M.$$ For example, if
$r=0$, then the $M$ functions $p_i(x)$ for $i=M+1,\ldots,2M$
generate an M-dimensional non-degenerate subspace of $W$ for the
form $r\Phi_1+s\Phi_2$. If $s=0$, then $p_i(x)$ for $i=1,\ldots,M$
generates an M-dimensional non-degenerate subspace for
$r\Phi_1+s\Phi_2$. And if neither $r$ nor $s$ is zero, either the
$M$ functions $p_i(x)$ for $i=M+1,\ldots,2M$ or for $i=1,\ldots,M$
generate an M-dimensional non-degenerate subspace.

Hence, any pencil of $\Phi_1$ and $\Phi_2$ has rank $\geq F(2)$,
since $M\geq F(2)$. Then, by Proposition \ref{prop:appro}, it has
the weak approximation, since $K$ is totally imaginary. Therefore,
the set of $K$-rational points in the intersection of $\Phi_1$ and
$\Phi_2$ on a non-degenerate subspace of dimension $\geq M$ is
Zariski-dense in the variety defined by $\Phi_1$ and $\Phi_2$. Let
$L$ be the hyperplane at $\infty$ and $L'$ the hyperplane defined
by $h(-2P)$ in the above notation of $f'=lh^2$ in case I or II. By
Lemma \ref{lem:intersec}, the intersection of two forms $\Phi_1$
and $\Phi_2$ is not contained in the union of two hyperplanes
defined by $L$ and $L'$. Hence, by the density of $K$-rational
points, we can get a nontrivial common zero over $K$ which is a
common zero of two original quadratic equations which defines an
elliptic function $f$ such that $f'=lh^2$, for some elliptic
function $h$ such that $h(-2P)\neq 0$.

 Now we take
the projective closure $V$ over $K$ of the linear subspace of $
\bbp(H^0(E,\mathcal{L}(n(O))))$ generated by $f$ and the constant
function $1$ over $K$. Note that the linear space generated by $f$
and $1$ is a base-point free linear system on $E$ from the
construction. Then $V$ is isomorphic to the projective line
$\bbp^1(K)$. And the normalization $X \subseteq E^{n-1}/A_n$ of
its preimage  under the 2-1 map from $E^{n-1}/A_n$ to
$E^{n-1}/S_n$ meets the ramification divisor wherever the divisor
$f-\lambda$ for some $\lambda$ has a zero or a pole of even
multiplicity $\geq 2$, that is, wherever its derivative
$(f-\lambda)'=f'$ has a zero or a pole of odd order. And $X$ has
only two points which meet the ramification locus  at $-2P$ and
$O$ to odd contact order by Lemma \ref{lem:contact} below. Hence
by the Hurwitz formula, the normalization of $X$ in $E^{n-1}/A_n$
is a curve of genus 0 defined over $K$. By subtracting the
constant $\lambda_p=f(-2P)$ from $f$, the function $f-\lambda_p$
has one double zero at $-2P$ and other zeros of odd order, since
$f'$ has only one simple zero at $-2P$ and other zeros of even
order.
\end{proof}

\begin{lem}\label{lem:contact} Under the same notation as in the proof of
Proposition \ref{prop:An}, if an elliptic function $f$ has a zero
(or a pole) at a point $P$ of order $m$, the contact order of $f$
with the ramification locus of the double cover from $E^{n-1}/A_n$
onto $E^{n-1}/S_n$ at $P$ is $m-1$.
\end{lem}

\begin{proof}
 Suppose $f$ has a zero $\alpha$ corresponding to the zero $P$ of order $m$.
Let $f(z)=(z-\alpha)^m(a_0+a_1(z-\alpha)+\cdot\cdot\cdot +$ higher
terms in $(z-\alpha))$, where $a_0\neq 0$.

 Note that the
ramification locus under the quotient map from $E^{n-1}$ to
$E^{n-1}/S_n$ is the zero locus of $\prod\limits_{i<j}(z_i-z_j)$,
where $z_i$ are zeros of $f-\lambda$,  for a parameter $\lambda$,
that is, the quotient map is ramified whenever $f-\lambda$ has a
double zero. By considering the ramification index, since the
degree of the map from $E^{n-1}/A_n$ onto $E^{n-1}/S_n$ is 2, the
ramification locus under the double cover from $E^{n-1}/A_n$ onto
$E^{n-1}/S_n$ is the zero locus of the discriminant of
$f-\lambda$, that is,
$$\prod\limits_{i<j}(z_i-z_j)^2,$$ where $z_i$ are zeros of
$f-\lambda$. If we write the discriminant of $f-\lambda$ in terms
of $\lambda$  with small $|\lambda|$, then its degree with respect
to $\lambda$ is the contact order of $f$ at $P$. We may assume
that $\alpha=0$ by translation. Hence we have
$$f-\lambda=0 \Leftrightarrow
z^m(a_0+a_1z+a_2z+\cdot\cdot\cdot+\mbox{ higher terms in } z
)-\lambda =0.$$ By Hensel's Lemma (\cite{sil86}, Chapter IV, Lemma
1.2) on $\bbc[[\lambda^{\frac{1}{m}}]]$,
%$\bbc\llbracket\lambda^{\frac{1}{m}}\rrbracket$,
all zeros of $z^m(a_0+a_1z+a_2z+\cdot\cdot\cdot+\mbox{ higher
terms in } z )-\lambda$ are
$$z_i=\left(\frac{\lambda}{a}\right)^{\frac{1}{m}} \zeta _m^i+A_i(\lambda), \mbox{for } 0\leq i\leq m-1,$$
where  $\zeta_m$ is a primitive $m$th root of unity, and
$A_i(\lambda)$ is a convergent Puiseux series in $\lambda$, that
is, a convergent power series in $\lambda^{\frac{1}{m}}$. Hence
the degree of  the discriminant of $f-\lambda$ with respect to
$\lambda$ is $\frac{1}{m} \cdot{ m\choose 2}\cdot 2=m-1$, which is
the contact order at $\alpha$ with the ramification locus. For a
pole, we proceed similarly, replacing $f$ by $1/f$.
\end{proof}

Next, we examine the Galois theory of the fixed fields
$\overline{K}^{\sigma}$ for automorphisms $\sigma\in
Gal(\overline{K}/K)$. We give some definitions.

\begin{defnn}\label{def:real} A field $F$ is $($formally$)$ real, if $-1$ is not a sum of squares
in $F$. A real field $F$ is real closed, if no algebraic extension
of $F$ is real. \end{defnn}

\begin{lem}\label{lem:brauer} Let $K$ be a number field. Then for any $\sigma\in
Gal(\overline{K}/K)$,
$$Gal(\overline{K}/\overline{K}^{\sigma}) \cong
\prod\limits_{p\in S}\bbz_p \mbox{\hspace{.2 in} or\hspace{.2in}}
\bbz/2\bbz,$$ where $S$ is a set of prime integers. In particular,
if $K$ is totally imaginary,
$Gal(\overline{K}/\overline{K}^{\sigma})$ has no torsion element,
hence, the Brauer group $Br(\overline{K}^{\sigma})$ of the fixed
field under $\sigma$ is trivial.
\end{lem}

\begin{proof} Let $\sigma\in Gal(\overline{K}/K)$. $Gal(\overline{K}/\overline{K}^{\sigma})$
is isomorphic to the closure of the subgroup generated by $\sigma$
in the sense of the Krull topology by (\cite{morandi}, Theorem
17.7). Hence,
$$Gal(\overline{K}/\overline{K}^{\sigma}) \cong
\prod\limits_{p\in S}\bbz_p\times \prod\limits_{p\in
T}\bbz/p^{m_p}\bbz \cong \prod\limits_{p\in S}\langle
\sigma_p\rangle \times \prod\limits_{p\in T}\langle
\tau_p\rangle,$$ where $S$ and $T$ are disjoint sets of primes,
$m_p$ are positive integers, and $\tau_p$ has a finite order
$p^{m_p}$. But since any element in
$Gal(\overline{K}/\overline{K}^{\sigma})$ has the order $1, 2$ or
$\infty$ by Artin-Schreier Theorem (\cite{issacs}, Theorem
(25.1)), the torsion part of
$Gal(\overline{K}/\overline{K}^{\sigma})$ is trivial or
$\bbz/2\bbz$. Moreover, if there are $q\in T$ and $p\in S\cup
(T-\{q\})$, then, $\tau_q$ is an involution and its fixed field
$\overline{K}^{\tau_q}$ is a real closed field by (\cite{issacs},
Theorem (25.13)). Also $\sigma_p^{-1}\tau_q\sigma_p=\tau_q$, so
$\sigma_p$ induces a nontrivial automorphism of
$\overline{K}^{\tau_q}$. This contradicts the uniqueness of an
isomorphism between two real closed fields (\cite{l93}, XI, \S2,
Theorem 2.9, pp.455). Therefore,
$$Gal(\overline{K}/\overline{K}^{\sigma})\cong \prod\limits_{p\in S}\bbz_p\mbox{\hspace{.2 in}or \hspace{.2 in}}
\bbz/ 2\bbz.$$

 On the other hand, if $Gal(\overline{K}/\overline{K}^{\sigma})$  is isomorphic to $
\bbz/ 2\bbz$ generated by $\tau$, then
$[\overline{K}:\overline{K}^{\tau}]=2$, so $\overline{K}^{\tau}$
is real-closed by (\cite{issacs}, Theorem (25.13)), so it has a
real embedding by (\cite{issacs}, Theorem (25.18)). Hence if
$K\subseteq\overline{K}^{\tau}$  is
 totally imaginary, then
  $Gal(\overline{K}/\overline{K}^{\sigma})$ is isomorphic to  $\prod\limits_{p\in S}\bbz_p$.
Then, since $\overline{K}^*$ is a divisible topological
$\prod\limits_{p\in S}\bbz_p$-group, $H^2(\prod\limits_{p\in
S}\bbz_p, \overline{K}^*)$  is trivial by (\cite{nsw}, Chapter 1,
\S 6. Proposition 1.6.13.(ii)). Therefore,
$Br(\overline{K}^{\sigma})=0$.
\end{proof}

\begin{lem}\label{lem:conic} Let $K$ be a totally imaginary number field. Then for any
$\sigma\in Gal(\overline{K}/K)$, a conic curve $X$ defined over
$K$ has a $\overline{K}^{\sigma}$-rational point.
\end{lem}

\begin{proof} Let $\sigma\in Gal(\overline{K}/K)$. Since a
conic can be identified with an element of
$Br(\overline{K}^{\sigma})$ as a Severi-Brauer variety of
dimension 1, and $Br(\overline{K}^{\sigma})=0$ by  Lemma
\ref{lem:brauer}, a conic is isomorphic to $\bbp^1$ over
$\overline{K}^{\sigma}$. Equivalently, it has  a
$\overline{K}^{\sigma}$-rational point.
\end{proof}

Let $f\in K(t_1,\ldots,t_m)[X_1,\ldots,X_n]$ be a polynomial with
coefficients in the quotient field $K(t_1,\ldots,t_m)$ of
$K[t_1,\ldots,t_m]$ which is irreducible over $K(t_1,\ldots,t_m)$.
We define $$H_K(f)=\{(a_1,\ldots,a_m)\in K^m\mid
f(a_1,\ldots,a_m,X_1,\ldots,X_n) \mbox{ is irreducible over }
K\}$$ to be the Hilbert set of $f$ over $K$. We need the following
lemma.

\begin{lem}\label{lem:hilbert}
Let $L$ be a finite separable extension of $K$ and let $f\in
L(t_1,\ldots,t_m)[X_1,\ldots,X_n]$ is an irreducible polynomial
over the quotient field $L(t_1,\ldots,t_m)$. Then, there exists a
polynomial $p\in K[t_1,\ldots,t_m,X_1,\ldots,X_n]$ such that $p$
is irreducible over $K(t_1,\ldots,t_m)$ and $H_K(p) \subseteq
H_L(f)$.
\end{lem}

\begin{proof} For a given irreducible polynomial
$f \in L(t_1,\ldots,t_m)[X_1,\ldots,X_n]$, by (\cite{jar}, Ch.11,
Lemma 11.6),
 there is an irreducible polynomial $q$  $\in
K(t_1,\ldots,t_m)[X_1,\ldots,X_n]$ such that $H_K(q)\subseteq
H_L(f)$. By (\cite{jar}, Ch.11, Lemma 11.1), there is an
irreducible polynomial $p \in K[t_1,\ldots,t_m,X_1,\ldots,X_n]$
which is irreducible over $K(t_1,\ldots,t_m)$ such that $H_K(p)
\subseteq H_K(q)$. Hence the Hilbert set $H_L(f)$ of $f$ over $L$
contains the Hilbert set $H_K(p)$ of $p$ over $K$.
\end{proof}

Let $G$ be a finite group and $\Lambda$ an $n$-dimensional
$G$-representation. Then, $G$ acts on $E\otimes \Lambda$ through
its action on $\Lambda$. Define $E\otimes\Lambda$ to be the
abelian variety representing the functor $S \mapsto
E(S)\otimes_{\bbz}\Lambda$, where $S$ is any scheme over the
ground field and $E(S)$ is the functor of points associated to
$E$. Then, as an abelian variety, $E\otimes\Lambda$ is just $E^n$,
since the action of $G$ on $E\otimes\Lambda$ is only though
$\Lambda$. With this background, we prove the following
proposition.

\begin{prop}\label{prop:group}
Let $K$ be a totally imaginary number field and $\sigma\in
Gal(\overline{K}/K)$. Let $G$ be a nontrivial finite group and
$\Lambda$ an $n$-dimensional  integral  $G$-representation
 for a positive integer $n$. Then $G$ acts on $E^n\cong E\otimes \Lambda$ through $\Lambda$.
 Suppose that there is a curve $X$ of genus 0 in $E^n/G$ over $K$.
Suppose the preimage of $X$ under the quotient map by $G$ is
decomposed into $k$ irreducible curves $C_1,\ldots, C_k$ such that
each $C_i\rightarrow X$ is the Galois covering with $Gal(C_i/X)
\leq G$, then $X$ cannot be decomposed completely, i.e. $k<|G|$,
and $Gal(C_i/X)$ are conjugate to each other in $G$. And
 for an irreducible component  $C\subseteq E^n$ in
  the preimage of $X$,
 there exist a finite extension $F$ of $K$ and an infinite sequence
 $\{L_i/F\}_{i=1}^{\infty}$ of linearly disjoint
finite Galois extensions of $F$ such that $F\subseteq
\overline{K}^{\sigma}$ and  $Gal(L_i/F)$ is naturally isomorphic
to $Gal(C/X)$ as a subgroup  of $G$. And for each $i$, there is a
submodule $M_i$ of $E(L_i)\otimes \bbq$ isomorphic to
$\Lambda\otimes \bbq$ as a $Gal(L_i/F)$-module via the inclusion
$Gal(L_i/F) \hookrightarrow G$.

In particular, if  $K$ is an arbitrary number field and $X$ is
isomorphic to $\bbp^1$ over $K$, then this holds with $F=K$. And
if the preimage of $X$ in $E^n$ is irreducible, then each
$Gal(L_i/F)$ is isomorphic to $G$ itself.
\end{prop}

\begin{proof} Let $\sigma\in Gal(\overline{K}/K)$. If the curve $X$ of genus 0  has a $K$-rational
 point, then $X\cong\bbp^1$ over $K$. If not, it is isomorphic to a conic curve. Then, since $K$ is
 totally imaginary, by Lemma \ref{lem:conic}, for every $\sigma\in Gal(\overline{K}/K)$, $X$ has a
 $\overline{K}^{\sigma}$-rational point.  Choose a point
 of $X$ over $\overline{K}^{\sigma}$ and let $F$ be the field of definition of this point. Then
 $F\subseteq \overline{K}^{\sigma}$, $F$ is a finite extension of $K$, and $X$ is isomorphic to $\bbp^1$ over $F$.
 Now we consider $X\subseteq E^n/G$ as $\bbp^1$ over $F$. Note that if $X\cong \bbp^1$ over $K$, then we can
 take $F=K$.

First, suppose that the preimage of the curve $X$ in $E^n$ under
the quotient map by $G$ is an irreducible curve $C$ with the
function field $F(C)$. Then the restricted quotient map
$\phi:C\rightarrow X$ by $G$ realizes $F(C)$ as a Galois extension
of the function field $F(x)$  of $X(F)\cong \bbp^1_F$ with the
Galois group isomorphic to $G$. By the theorem of the primitive
element, there exists $t\in F(C)$ such that $F(C)=F(x,t)$ and
$$g_mt^{m}+g_{m-1}t^{m-1}+\cdot\cdot\cdot+g_1t+g_0=0,$$
where $g_i$ are polynomials in $F[x]$.  Choose a minimal
polynomial of $t$ over $F(x)$ and clear its denominators so that
we let $f(x,y)$ be a minimal polynomial of $t$ in $F[x,y]$. Then,
$f$ is absolutely irreducible over $F$ so it is irreducible over
$F(x)$.

By (\cite{sil}, Lemma), the set $\bigcup\limits_{[L:F]\leq k}
E(L)_{tor}$ is a finite set, where the union runs over all finite
extensions $L$ of $F$ whose degree over $F$ is $\leq k$, where
$k=|G|$. Let $L'$ be a finite field extension of $F$ over which
all points of $\bigcup\limits_{[L:F]\leq k} E(L)_{tor}$ are
defined. Applying Lemma \ref{lem:hilbert} and (\cite{jar}, Lemma
12.12) to $f$ over $L'$, we can choose $x_1\in H_F(f)\cap K$ such
that the specialization $x \mapsto x_1$ preserves the Galois group
$G$ and there is a point $Q_1$ of $C\subseteq E^n\cong E\otimes
\Lambda $ in the preimage $\phi^{-1}((1:x_1))$ of $(1:x_1)\in
\bbp^1(F)\cong X$ under $\phi$ is defined over a finite Galois
extension $L_1$ of $F$ with $Gal(L_1/F)\cong G$, that is, the
preimage of $(1:x_1)$ under $\phi$ consists of a single point Spec
$L_1$. Let $\Lambda^*$ be the dual of $\Lambda$ with the action of
$G$. Then the morphism from Spec $L_1$  to $E\otimes \Lambda$
induces a $\bbz[G]$-linear map $g :\Lambda^{*} \rightarrow E(L_1)$
given by $\lambda^*\mapsto \sum\limits_j\lambda^*(\lambda_j)P_j$,
where $Q_1=\sum\limits_j P_i\otimes \lambda_i\in E\otimes
\Lambda$.

Since $f(x_1,y)$ is irreducible over $L'$, two extensions $L_1$
and $L'$ are linearly disjoint over $F$. So for $\lambda^*\in
\Lambda^*$, $g(\lambda^*)\in E(L_1)$ is a non-torsion point. So if
we let $M_1\subseteq E(L_1)\otimes \bbc$ be the submodule
generated by the points of $E(L_1)$ in the image of $\Lambda^*$
under the given map $g$ in the above, it is a submodule of
$E(L_1)\otimes \bbq$ isomorphic to $\Lambda^*\otimes \bbq$ as a
$Gal(L_1/F)$-module via the natural isomorphism $Gal(L_1/F)\cong
G$. Since $\Lambda$ is a finite dimensional integral
representation, it is isomorphic to its dual $\Lambda^*$ as
$G$-representations. So $M_1$ is isomorphic to $\Lambda\otimes
\bbq$ as a $Gal(L_1/F)$-module.

Suppose the preimage of $X$ is decomposed into a union of
irreducible curves $C_1\cup C_2\cup\cdot\cdot\cdot\cup C_k$. Then,
$G$ acts transitively on the set of $k$ curves and each
$Gal(C_i/K)$ can be identified with the stabilizer of $C_i$ in $G$
so $Gal(C_i/K)$ are conjugate to each other. So if $k=|G|$, then
this implies that $C_i\cong \bbp^1$ in $E^n$ which is impossible,
because no abelian variety  contains $\bbp^1$ as a subvariety. So
$k <|G|$.  Let $C$ be one of irreducible components $C_i$.
Applying the same argument with the quotient map from the fixed
component $C$ to $X$, we get a Galois extension $L_1$ of $F$ with
the Galois group $Gal(L_1/F)$ which is isomorphic to the
stabilizer of $C$ in $G$ which is $Gal(C/X)\leq G$ and a
$Gal(L_1/F)$-submodules $M_1$ of $E(L)\otimes \bbq$ generated by
$n$ non-torsion points of $E(L_1)$ and
 it is isomorphic to $\Lambda\otimes \bbq$ as a
$Gal(L_1/F)$-module via the natural inclusion
$Gal(L_1/F)\hookrightarrow G$.

Inductively, suppose we have found  linearly disjoint finite
Galois extensions $L_1,L_2,\ldots,L_{k}$  of $F$ and  for each
$i=1,2,\ldots,k$, there is a submodule $M_i$ of $E(L_i)\otimes
\bbq$ isomorphic to $\Lambda\otimes\bbq$ as a $Gal(L_i/F)$-module
via the natural inclusion $Gal(L_1/F)\hookrightarrow G$. By
applying Lemma \ref{lem:hilbert} and (\cite{jar}, Lemma 12.12) to
$f$ over the composite field $L'L_1L_2\cdot\cdot\cdot L_k$, there
is a point $x_{k+1}\in X(F)$ such that the specialization $x
\mapsto x_{k+1}$ preserves the Galois group $G$ and a point in the
preimage of $x_{k+1}$ in $C$ is defined over a Galois extension
$L_{k+1}$ of $F$ which is linearly disjoint from
$L'L_1L_2\cdot\cdot\cdot L_k$ and has the Galois group isomorphic
to a subgroup of $G$. Then similarly, we get  a
$Gal(L_{k+1}/F)$-submodule $M_{k+1}$ generated by $n$  non-torsion
points of $E(L_{k+1})$ isomorphic to $\Lambda\otimes \bbq$ via
$Gal(L_{k+1}/F)\hookrightarrow G$. This completes the proof.
\end{proof}

\begin{cor}\label{cor:group}
Let $K$ be a totally imaginary number field and $E/K$ an elliptic
curve over $K$ with a $K$-rational point such that $2P\neq O$ and
$3P\neq O$. Let $\Lambda$ be the $(n-1)$-dimensional irreducible
quotient representation space of the natural permutation
representation of the alternating group $A_n$ by the trivial
representation. Let $\sigma\in Gal(\overline{K}/K)$. Then for some
even integer $n$,
 there exist a finite extension $F\subseteq \overline{K}^{\sigma}$ over $K$ and an
infinite  sequence $\{L_i/F\}_{i=1}^{\infty}$ of linearly disjoint
finite Galois extensions of $F$ such that $Gal(L_i/F)$ acts
transitively on $\{1,2,\ldots,n\}$ as a subgroup of $A_n$. And for
each Galois extension $L_i$ of $K$, there is a submodule $M_i$ of
$E(L_i)\otimes \bbq$ isomorphic to the $(n-1)$-dimensional
irreducible quotient representation space $\Lambda\otimes\bbq$ as
a $Gal(L_i/F)$-module via the natural inclusion
$Gal(L_i/F)\hookrightarrow A_n$.
\end{cor}

\begin{proof} By Proposition \ref{prop:An}, there is a curve
$X$ of genus 0 defined over $K$ in $E^{n-1}/A_n$, for some even
integer $n$. So by Proposition \ref{prop:group}, there exist such
an infinite sequence of Galois extensions $L_i$ and submodules
$M_i$ of $E(L_i)\otimes \bbq$ isomorphic to the
$(n-1)$-dimensional irreducible quotient representation space
$\Lambda\otimes\bbq$ of $A_n$ as a $Gal(L_i/F)$-module via the
natural inclusion $Gal(L_i/F)\hookrightarrow A_n$.
 And by Proposition \ref{prop:An}, Proposition \ref{prop:group}, and Lemma \ref{lem:tran},
for each $Gal(L_i/F)$ as a subgroup of $A_n$, there is a subgroup
$H_i\leq S_n$ such that $H_i\cap A_n\cong Gal(L_i/F)$ and it acts
transitively on $\{1,2,\ldots,n\}$. Moreover, the image of $X$
given by Proposition \ref{prop:An} under the 2-to-1 map from
$E^{n-1}/A_n$ onto $E^{n-1}/S_n$  has a divisor which decomposes
into one divisor of ramification degree 2 and other divisors of
 odd degree. So by Lemma \ref{lem:(12)}, $H_i$ contains a
transposition. Therefore, by Lemma \ref{lem:semidirect},
 $Gal(L_i/F)$ also acts transitively on $\{1,2,\ldots,n\}$.
\end{proof}

\begin{lem}\label{lem:indep}
Let $E/K$ be an elliptic curve over a number field $K$. Let $d$ be
a positive integer $\geq 2$. Suppose $\{L_i/K\}^{\infty}_{i=1}$ is
an infinite sequence of linearly disjoint finite Galois extensions
of $K$ whose degrees $[L_i:K]$ are $\leq d$ and
$\{P_i\}^{\infty}_{i=1}$ is an infinite sequence of points in
$E(\overline{K})$ such that for each $i$, $P_i\in E(L_i)$ but
$P_i\notin E(K)$. Then, there is an integer $N$ such that
$\{P_i\}_{i\geq N}$ is a sequence of linearly independent
non-torsion points of $E$.
\end{lem}

\begin{proof} By (\cite{sil}, Lemma), the set $S=\bigcup\limits_{[L:K]\leq d}
E(L)_{tor}$ is a finite set, where the union runs all over finite
extensions $L$ of $K$ whose degree over $K$ is $\leq d$. So there
is a finite extension $F$ of $K$ over which all points of $S$ are
defined and there is an integer $n$ such that $nP=O$, for all
$P\in S$. Let $n$ be such a fixed integer and let $T$ be the set
of all points $P$ of $E(\overline{K})$ such that $n P\in E(K)$.
Then, since $E(K)$ is finitely generated by the Mordell-Weil
Theorem (\cite{sil86}, Chapter VIII), there is a finite extension
$F'$ of $K$ over which all points of $T$ are defined. Then, all
but finitely many fields $L_i$ in the given sequence
$\{L_i/K\}_{i=1}^{\infty}$ are linearly disjoint from $F$ and $F'$
over $K$. This implies that there is an integer $N$ such that
points $P_i$ for all $i\geq N$ are non-torsion points in $E(L_i)$.
And by linear disjointness of fields $L_i$, $F$ and $F'$ over $K$,
we have that for all $i\geq N$,
$$E(L_i)\cap S \subseteq E(K)_{tor} \mbox{~~and ~~} E(L_i)\cap T\subseteq E(K).$$
Note that since each $P_i\notin E(K)$, we have that for any
integer $m\geq N$ and for each $i$ such that $N\leq i\leq m$,
there is an automorphism $\tau_i\in Gal(\overline{K}/K)$ such that
$\tau_i|_{L_j}=id_{L_j}$ for all $N\leq j\neq i\leq m$, but
$\tau_i(P_i)\neq P_i$. Moreover, we may choose such a $\tau_i$
that $\tau_i(P_i)-P_i$ is not a torsion point. In fact, otherwise,
for every restriction $\tau_i|_{L_i}\in Gal(L_i/K)$ of $\tau_i$,
$\tau_i|_{L_i}(P_i)-P_i$ is a torsion point in $E(L_i)$. Hence,
$\tau_i|_{L_i}(P_i)-P_i \in E(L_i)\cap S \subseteq E(K)_{tor}.$
Then, $n(\tau_i|_{L_i}(P_i)-P_i)=O$ so $\tau_i|_{L_i}(nP_i)=nP_i$
for all $\tau_i|_{L_i}\in Gal(L_i/K)$. This implies $nP_i\in E(K)$
so $P_i\in T\cap E(L_i)\subseteq E(K)$ which contradicts the
assumption that $P_i\notin E(K)$. Hence, we conclude that for each
$i$ such that $N\leq i\leq m$, there is an automorphism $\tau_i\in
Gal(\overline{K}/K)$ such that $\tau_i|_{L_j}=id_{L_j}$ for all
$N\leq j\neq i\leq m$, but $\tau_i(P_i)-P_i$ is a non-trivial and
non-torsion point of $E$.

Let $m\geq N$ be a given positive integer. Suppose that for some
integers $a_i$,
$$a_NP_{N}+a_{N+1}P_{N+1}+\cdots+a_mP_{m}=0.$$ By the claim
 above, for each $i=N, N+1, \ldots, m$, there is an
automorphism $\tau_i\in Gal(\overline{K}/K)$ such that
$\tau_i|_{L_{j}}=id_{L_{j}}$ for all $1\leq j\neq i \leq m$ but
$\tau_i(P_i)-P_i$ is a non-trivial and non-torsion point of $E$.
Now we apply such $\tau_i$ to get
$$a_NP_{N}+a_{N+1}P_{N+1}+\cdots+a_{i-1}P_{i-1}+a_i\tau_i(P_{i})+a_{i+1}P_{i+1}+\cdots+a_mP_{i_m}=0.$$
So by subtracting, we get $a_i(P_{i}-\tau_i(P_{i}))=0$, which
implies $a_i=0$. Hence any non-torsion points in $\{P_{i}\}_{i\geq
N}$ are linearly independent.
\end{proof}

\begin{thm}\label{thm:totally}  Let $K$ be a totally imaginary number field. Suppose  $E/K$ has a
$K$-rational point $P$ such that $2P\neq O$ and $3P\neq O$. Then
for each $\sigma \in Gal(\overline{K}/K)$,
$E(\overline{K}^{\sigma})$ has infinite rank.
\end{thm}

\begin{proof}  Let $\sigma\in Gal(\overline{K}/K)$.
By Proposition \ref{prop:An} and Corollary \ref{cor:group}, there
are a finite extension $F\subseteq \overline{K}^{\sigma}$ over $K$
and an infinite sequence $\{L_i/F\}_{i=1}^{\infty}$ of linearly
disjoint finite Galois extensions of $F$ such that the Galois
group $Gal(L_i/F)$ acts transitively on $\{1,2,\ldots,n\}$ as a
subgroup of $A_n$ for some even integer $n$. And for each $i$,
there is a $Gal(L_i/F)$-submodule of $E(L_i)\otimes\bbq$ which is
isomorphic to the restriction of the natural ($n-1$)-dimensional
quotient of the permutation representation  of $A_n$ to
$Gal(L_i/F)$.

 Let $\sigma _i = \sigma|_{L_i}$ be the restriction of $\sigma$ to $L_i$.
Then  since $F\subseteq \overline{K}^{\sigma}$,
$\sigma_i|_F=id_F$. Therefore, $\sigma_i\in Gal(L_i/F)\leq A_n$.
Let $E(L_i^{\sigma _i})$ be the group of fixed points of $E(L_i)$
under $\sigma_i$. Then obviously, $E(L_i^{\sigma _i})\subseteq
E(\overline{K}^{\sigma})$. Since $n$ is even, each $M_i$ of
$E(L_i)\otimes\bbq$
 has a fixed element $v_i$ under $\sigma_i$
by Lemma \ref{lem:cycle} and Lemma \ref{lem:inva}.

 Note that each $v_i$ is not a torsion point and not defined
over $F$. In fact, if $v_i$ is defined over $F$, then $v_i$ is
fixed under every element in $Gal(L_i/F)$. But since $Gal(L_i/F)$
acts transitively on $\{1,2,\ldots,n\}$,  by Lemma
\ref{lem:invari}, there is no fixed vector of the restriction to
$Gal(L_i/F)$ of the $(n-1)$-dimensional quotient of the
permutation representation of $A_n\leq S_n$. Then, by Lemma
\ref{lem:indep}, there is an integer $N$ such that $\{v_i\}_{i\geq
N}$ are linearly independent.

 Since $v_i\in E(L_i^{\sigma _i})\otimes\bbq$ for each $i$, the module generated by
$\{v_i\}_{i=1}^{\infty}$ over $\bbq$ is a submodule of
$E\left(\prod \limits _{i=1}^{\infty} L_i^{\sigma
_i}\right)\otimes\bbq$. Hence
$$E(\overline{K}^{\sigma})\otimes\bbq\supseteq E\left(\prod
\limits _{i=1}^{\infty} L_i(\sigma _i)\right)\otimes\bbq \supseteq
\{v_i\}_{i=1}^{\infty}\supseteq \{v_i\}_{i\geq N}$$ is infinite
dimensional.
\end{proof}

\section{Infinite rank over the fixed fields under complex conjugation automorphisms}

The only difficulty in proving the rank of
$E(\overline{K}^{\sigma})$ is infinite is that
$\overline{K}^{\sigma}$ may have a real embedding. Now we consider
complex conjugation automorphisms of $\overline{K}$ and
 prove that without hypothesis on rational points of elliptic curves and the ground field,
 the rank of an elliptic curve over the fixed field under every complex conjugation automorphism is
 infinite.

\begin{defnn} A field $F$ is
called an ordered field with the positive set $P$, if
$F=P\bigsqcup \{0\}\bigsqcup -P$, a disjoint union, where $P$ is a
subset of $F$ closed under addition and multiplication.
\end{defnn}

Now we prove the following two lemmas by using the relation
between real fields (refer Definition \ref{def:real}) and ordered
fields.

\begin{lem}\label{lem:real1} If a field $F$ is ordered (or real)  and
algebraic over $\bbq$, then $F$ has a real embedding $\theta$,
that is, $\theta(F) \subseteq \bbr\cap \overline{F}$.
\end{lem}

\begin{proof} By (\cite{issacs}, Chapter
25, Corollary (25.22), pp.411), a field $F$ is ordered if and only
if it is real. Hence, $F$ is real. And since $F$ is real and
algebraic over $\bbq$, by (\cite{issacs}, Theorem (25.18),
pp.~410), there exists an isomorphism from $F$ into $\bbr\cap
\overline{F}$.
\end{proof}

 We give some
equivalent statements of complex conjugation automorphisms.

\begin{lem}\label{lem:real2} The following statements are
equivalent: for an automorphism $\sigma\in Gal(\overline{K}/K)$, \\
(1) $\overline{K}^{\sigma}$ has a real embedding $\theta$, that
is,
$\theta(\overline{K}^{\sigma}) \subseteq \bbr\cap\overline{K}$\\
(2) $\sigma$ is a complex conjugation automorphism, that is, the
order of $\sigma$ in $Gal(\overline{K}/K)$ is 2.\\
(3) $\overline{K}^{\sigma} \cong\bbr\cap\overline{K}$
\end{lem}

\begin{proof} Suppose (1). Then,
$$\langle \sigma \rangle \cong
Gal(\overline{K}/\overline{K}^{\sigma})\cong
Gal(\overline{K}/\theta(\overline{K}^{\sigma})) \unrhd
Gal(\overline{K}/\bbr\cap \overline{K}) \cong \bbz/2\bbz,$$ since
$[\overline{K}:\bbr\cap \overline{K}]=2$. Hence,
$Gal(\overline{K}/\overline{K}^{\sigma})$ has a torsion subgroup
of order 2. Then, by Lemma \ref{lem:brauer},
$Gal(\overline{K}/\overline{K}^{\sigma})$ itself is isomorphic to
$\bbz/2\bbz$. Since $\sigma$ is not trivial, we have
$Gal(\overline{K}/\overline{K}^{\sigma}) \cong
Gal(\overline{K}/\bbr\cap \overline{K})$, hence, the order of
$\sigma$ is $2$, which implies (2). And $\overline{K}^{\sigma}
\cong \bbr\cap \overline{K}$, which implies (3).

 Now we suppose (3). Then, the order of
$\sigma$ equals the degree $[\overline{K}:\overline{K}^{\sigma}]$
which is equal to $[\overline{K}:\bbr\cap\overline{K})]=2$, and
this implies (2).

Suppose (2). Then, $[\overline{K}:\overline{K}^{\sigma}]=2$.  By
by (\cite{issacs}, Theorem (25.13)), $\overline{K}^{\sigma}$ is
real closed. Then, it is real and algebraic over $\bbq$. So by
Lemma \ref{lem:real1}, it has a real embedding. This implies (1).
\end{proof}

 The following lemma gives the density of the
Hilbert sets over a number field $K$  with respect to any real
embeddings of $K$ into $\bbr$.

\begin{lem}
\label{lem:dense} Let $K$ be a number field and
$\tau_1,\ldots,\tau_m$ be a family of real embeddings of $K$. For
$i=1,2,\ldots,k$, let $f_i(x,y)\in K[x,y]$ be irreducible
polynomials over $K(x)$. Let $H_K(f_i)$
%$= \{\alpha\in K\mid$ $ f_i(\alpha, y)
%\in K[y]$ is irreducible over $K\}$
be the Hilbert  set of $f_i$ over $K$. Then
$$\left(\bigcap\limits _{i=1}^k
H_K(f_i)\right)~\cap~\left(\bigcap\limits
_{j=1}^m\tau_j^{-1}(I)\right)\neq~~ \emptyset,$$ for any open
interval $I$ in $\bbr$.
\end{lem}

\begin{proof} This is a special case of (\cite{geyer},
Lemma 3.4).
\end{proof}

\begin{thm}
\label{thm:rank} Let $K$ be a number field, $K_{ab}$ the maximal
abelian extension of $K$,  and E/K an elliptic curve over $K$.
Then, for any complex conjugation automorphism $\sigma\in
Gal(\overline{K}/K)$, $E((K_{ab})^{\sigma})$ has infinite rank.
Hence, $E(\overline{K}^{\sigma})$ has infinite rank.
\end{thm}

\begin{proof}
For  a complex conjugation automorphism $\sigma\in
Gal(\overline{K}/K)$, there exists a real embedding $\theta$ such
that $\theta(\overline{K}^{\sigma})\subseteq \bbr\cap
\overline{K}$, by Lemma \ref{lem:real2}. Note that
$\sigma(\sqrt{-1})=-\sqrt{-1}$, since otherwise
$\sigma(\sqrt{-1})=\sqrt{-1}$ and then, $\sqrt{-1}\in
\overline{K}^\sigma$ and $0<(\theta(\sqrt{-1}))^2=\theta(-1)=-1$
which is a contradiction. And for any element $\alpha\in K$ such
that $\theta(\alpha)> 0$, $\sigma(\sqrt{\alpha})=\sqrt{\alpha}$,
since otherwise, $\sigma(\sqrt{\alpha})=-\sqrt{\alpha}$, hence
$\sigma(\sqrt{-\alpha})=\sqrt{-\alpha}$, then, $ \sqrt{-\alpha}\in
\overline{K}^\sigma$ and
$0<(\theta(\sqrt{-\alpha}))^2=-\theta(\alpha)<0$ which is a
contradiction.

Let's consider a Weierstrass equation of $E/K$, $y^2=x^3+ax+b$,
for $a,b\in K$. Then there exists $\alpha\in \bbr$ such that
$x^3+\theta(a)x+\theta(b)>0$ for all $x>\alpha$. Let
$I=(\alpha,\infty)$ be the open interval of all real numbers
$>\alpha$. Let $f(x,y)=y^2-(x^3+ax+b)$. Then $f$ is an absolutely
irreducible polynomial in $K[x,y]$, hence irreducible over $K(x)$.
Let $H_K(f)$ be the Hilbert set of $f$ over $K$. Note that the
restriction $\theta|_K$ of $\theta$ to $K$ is a real embedding of
$K$.

 By Lemma \ref{lem:dense}, there is an element $x_1\in  H_K(f)\cap \theta|_{K}^{-1}(I)$.
 Then $x_1^3+ax_1+b\in K $ and since $\theta(x_1)>\alpha$, $\theta(x_1^3+ax_1+b)$
is positive, hence $\sigma(\sqrt{x_1^3+ax_1+b})=
\sqrt{x_1^3+ax_1+b}$. Hence $\sigma$ fixes $\sqrt{x_1^3+ax_1+b}$.
Let $K_1=K(\sqrt{x_1^3+ax_1+b})$. Then since $f(x_1,y)$ is
irreducible in $K[y]$, $K_1$ is a quadratic extension of $K$ and
$K_1\subseteq\overline{K}(\sigma)$.

Inductively, suppose we have constructed linearly disjoint
quadratic extensions $K_1,$ $\ldots,K_{n-1}$ of $K$ such that
$K_i=K(\sqrt{x_i^3+ax_i+b})$ for $x_i\in H_K(f)\cap
\theta|_{K}^{-1}(I)$. Let $L_{n-1}=K_1\cdots K_{n-1}$ be the
composite field extension over $K$. By Lemma \ref{lem:dense}
again, there is $x_n\in H_{L_{n-1}}(f)\cap \theta|_K^{-1}(I)$. Let
$K_n$ $=K(\sqrt{x_n^3+ax_n+b})$. Then similarly, we can show that
$K_n$ is a quadratic extension of $K$ and
$K_n\subseteq\overline{K}(\sigma)$. Moreover,  $K_n$ is linearly
disjoint from all $K_1,K_2,\ldots,K_{n-1}$, since $x_n\in
H_{L_{n-1}}(f)$.

Hence we have obtained $\{x_i\}_{i=1}^{\infty}\subseteq K$ and an
infinite sequence $\{K_i/K\}_{i=1}^{\infty}$ of linearly disjoint
quadratic extensions of $K$ such that
$K_i=K(\sqrt{x_i^3+ax_i+b})$. For each $i$, let $P_i$ be a point
of $E(K_i)$ whose $x$-coordinate is $x_i$. Note that $P_i\notin
E(K)$, for each $i$. Hence, by Lemma \ref{lem:indep}, for some
$N$,  $\{P_i\}_{i\geq N}$ consists of  linearly independent
non-torsion points of $E(\overline{K})$. In particular, since
$K_i$ are abelian extensions of $K$ and are fixed under $\sigma$,
$P_i$ are points of $E((K_{ab})^{\sigma})$. Hence
$$E\left((K_{ab})^{\sigma}\right)\otimes \bbq \supseteq
E\left(\prod \limits _{i=N}^{\infty} K_i\right)\otimes \bbq
\supseteq \{P_i\otimes 1\}_{i\geq N}$$ is infinite dimensional.
And $E(\overline{K}^{\sigma})$ has infinite rank as well.
\end{proof}

\section{More general result : $E$ over arbitrary number fields with
a rational point which is neither 2-torsion nor 3-torsion }

In this section, we prove a more general result than the result of
Theorem \ref{thm:totally} for $E/K$ with a rational point $P$ such
that $2P\neq O$  and $3P\neq O$ without hypothesis on the ground
field $K$. To do so, we need the following lemma and proposition.

\begin{lem}\label{lem:real} For a number field $K$, let $\sigma\in
Gal(\overline{K}/K)$. If $\sigma$ does not fix any totally
imaginary finite extensions of $K$, then $\sigma$ is a complex
conjugation automorphism.
\end{lem}

\begin{proof}  Since $\overline{K}^{\sigma}$ is algebraic over $\bbq$, by Lemma \ref{lem:real1}
and Lemma \ref{lem:real2}, it is enough to show that
$\overline{K}^\sigma$ is ordered.

If $L$ is a finite extension of $K$ such that $L\subseteq
\overline{K}^\sigma$, then $L$ is not totally imaginary by the
assumption. Let $\tau_1,\ldots,\tau_r$ be all real embeddings of
$L$.

For $\alpha\in L^*$ ($=L-\{0\}$), if $\tau_i(\alpha)<0$ for all
$i=1,\ldots,r$, then, $L(\sqrt{\alpha})$ is totally imaginary
(otherwise, $L(\sqrt{\alpha})$ has a real embedding $\rho$ and
$\rho|_L=\tau_i$ for some $i$. But we have $0<
(\rho(\sqrt{\alpha}))^2=\rho(\alpha)=\tau_i(\alpha)$, which
contradicts $\tau_i(\alpha)<0$ for all $i$). Hence, $\sigma$ does
not fix $\sqrt{\alpha}$ by the assumption, so
$\sigma(\sqrt{\alpha})=-\sqrt{\alpha}$. This implies that for
$\beta\in L^*$,  if $\tau_i(\beta)>0$ for all $i=1,\ldots,r$, then
$\sigma(\sqrt{-\beta})=-\sqrt{-\beta}$ and since $\tau_i(-1)=-1<0$
for all $i$, $\sigma(\sqrt{-1})=-\sqrt{-1}$, hence
$\sigma(\sqrt{\beta})=\sqrt{\beta}$.

Therefore, there is a homomorphism $h: \prod\limits_{i=1}^r \{\pm
1\}\rightarrow \{\pm 1\}$ such that the action of $\sigma$ on
$\sqrt{\alpha}$ for $\alpha\in L^*$ depends only on the image of
the vector of signs of $\alpha$ under $h$. In other words, for
$\alpha\in L^*$, we let $f: L^*\rightarrow \prod\limits_{i=1}^r
\{\pm 1\}$ be a homomorphism defined by
$$f(\alpha)=(\mbox{sign}(\tau_1(\alpha)),\ldots,\mbox{sign}(\tau_r(\alpha))),$$
and $g: L^*\rightarrow \{\pm 1\}$ defined by $$g(\alpha)=\mbox{
the sign of }\frac{\sigma(\sqrt{\alpha})}{\sqrt{\alpha}}_,$$ so
$\sigma(\sqrt{\alpha})=g(\alpha)\sqrt{\alpha}$, then there exists
a homomorphism $h: \prod\limits_{i=1}^r \{\pm 1\} \rightarrow
\{\pm 1\}$ such that $h\circ f=g$.

Note that from the above explanation on totally positive or
totally negative elements of $L^*$, we get
$$(*)\hspace{1.5 in} h(-1,\ldots,-1)=-1 \mbox{\hspace{.2in} and
\hspace{.2in}} h(1,\ldots,1)=1.\hspace{1.5 in}$$ In particular,
there is always a vector consisting of $-1$ in all but one
coordinate and $1$ in the remaining coordinate which lies in the
kernel of $h$. In fact, there are $r$ vectors consisting of $-1$
in all but one coordinate and $1$ in the remaining coordinate:
$$v_1=(1,-1,\ldots,-1),v_2=(-1,1,-1,\ldots,-1),\ldots,v_r=(-1,\ldots,-1,1).$$
If all $r$ vectors map to $-1$ under $h$, then
$$(-1)^r=\prod\limits_{i=1}^r h(v_i)=h\left(\prod\limits_{i=1}^r
v_i\right)=h((-1)^{r-1},\ldots, (-1)^{r-1}).$$ But this
contradicts ($*$) by taking an even and odd integer $r$.
Therefore, at least one of $v_i$ must map to $1$, so it lies in
the kernel of $h$. Without loss of generality, we may assume that
$v_1$ maps to $1$ under $h$.

Hence, we can choose $\alpha\in L^*$ such that
$\sigma(\sqrt{\alpha})=\sqrt{\alpha}$ and $\tau_1(\alpha)>0$ but
$\tau_i(\alpha)<0$ for all $i=2,\ldots, r$ and let
$L'=L(\sqrt{\alpha})$. Then, $L'$ is fixed under $\sigma$ so $L'$
is not totally imaginary. Let $\rho$ be a real embedding of $L'$.
Then, since $\alpha$ is positive only with respect to $\tau_1$,
$$ 0<
(\rho(\sqrt{\alpha}))^2=\rho(\alpha)=\rho|_L(\alpha)=\tau_1(\alpha).$$
Hence, $$\rho(\sqrt{\alpha})=\pm \sqrt{\tau_1(\alpha)}.$$ This
shows that $L'$ has exactly two real embeddings $\rho_1$ and
$\rho_2$ such that
$$\rho_1(\sqrt{\alpha})=\sqrt{\tau_1(\alpha)}, \hspace{.15 in}
\rho_2(\sqrt{\alpha})=-\sqrt{\tau_1(\alpha)},\mbox{~ and ~}
\rho_{i}|_L=\tau_1, \mbox{ for } i=1,2.$$

We proceed the same argument on $L'$ with two real embeddings
$\rho_i$ as before and get a homomorphism $h':\{\pm 1\}\times
\{\pm 1\}\rightarrow \{\pm 1\}$ and $f': L'^*\rightarrow \{\pm
1\}\times \{\pm 1\}$ given by $f(\beta)=($sign$(\rho_1(\beta)),
$sign$(\rho_2(\beta)))$, $g':L'^*\rightarrow \{\pm 1\}$ given by
$g'(\beta)=\mbox{ the sign of
}\frac{\sigma(\sqrt{\beta})}{\sqrt{\beta}}$, such that $h'\circ
f'=g'$ on $L'^*$. Again, we have that
$$(**)\hspace{1.7 in} h'(-1,-1)=-1 \mbox{\hspace{.2in} and
\hspace{.2in}} h'(1,1)=1.\hspace{1.7 in}$$ This implies that $h'$
cannot send both $(-1,-1)$ and $(1,1)$ to the same value $1$ or
$-1$. So either $h'(-1,1)=-1$ and $h'(1,-1)=1$ or $h'(-1,1)=1$ and
$h'(1,-1)=-1$. Therefore, $h'$ is the projection onto either the
first factor or the second factor. Without loss generality, we
assume that $h'$ is the projection onto the first factor, that is,
$g'$ is defined by the sign of the first real embedding $\rho_1$.
Then,  if $\beta, \gamma\in L'^*$ such that $\sqrt{\beta},
\sqrt{\gamma}$ are fixed under $\sigma$, then $\rho_1(\beta)>0$
and $\rho_1(\gamma)>0$ so $\rho_1(\beta+\gamma)>0$. Hence,
$$g'(\beta+\gamma)=h'(f'(\beta+\gamma))=h'(1, a)=1, \mbox{ where } a= 1\mbox{ or } -1.$$ So
$\sigma(\sqrt{\beta+\gamma})=\sqrt{\beta+\gamma}.$ And obviously,
 $\sigma(\sqrt{\beta\gamma})=\sqrt{\beta\gamma}.$

 We have shown that the set of $\beta\in L'^*$ with
$\sigma(\sqrt{\beta})=\sqrt{\beta}$ is closed under addition and
multiplication. Therefore, for any two elements $a$ and $b\in
\overline{K}^\sigma-\{0\}$, by applying the preceding argument
with taking $L$ as a finite extension of $K$ generated by $a^2$
and $b^2$, the set of squares in $ \overline{K}^\sigma-\{0\}$ are
closed under addition and multiplication. Hence, if we let $S$ be
the set of squares in $ \overline{K}^\sigma-\{0\}$ and denote the
the set of non-squares in $ \overline{K}^\sigma-\{0\}$ by $-S$,
then $\overline{K}^\sigma= S \bigsqcup \{0\}\bigsqcup -S$, a
disjoint union. Hence, $\overline{K}^\sigma$ is ordered with the
positive set $S$. This completes the proof.
\end{proof}

\begin{prop}\label{prop:comptot} For a number field $K$, let $\sigma\in
Gal(\overline{K}/K)$. If $\sigma$ is not a complex conjugation
automorphism, then there is a totally imaginary finite extension
$L$ over $K$ such that $L \subseteq \overline{K}^\sigma$.
\end{prop}
\begin{proof} It follows
Lemma \ref{lem:real}.
\end{proof}

The following is our more general theorem without hypothesis on
the ground field or the given automorphisms.

\begin{thm}\label{thm:main} Let $K$ be a number field and $E/K$ an
elliptic curve over $K$ with a $K$-rational point $P$ such that
$2P\neq O$ and $3P\neq O$. Then, for each $\sigma\in
Gal(\overline{K}/K)$, the rank of $E(\overline{K}^\sigma)$ is
infinite.
\end{thm}

\begin{proof} Let $\sigma\in Gal(\overline{K}/K)$. If $\sigma$ is
a complex conjugation automorphism, then $E(\overline{K}^\sigma)$
has infinite rank by Theorem \ref{thm:rank}. If $\sigma$ is not a
complex conjugation automorphism, then by Proposition
\ref{prop:comptot}, there is a totally imaginary finite extension
$L$ of $K$ such that $L \subseteq \overline{K}^\sigma$. Hence
$\sigma\in Gal(\overline{K}/L)$. Now consider $E/L$ defined over
$L$ by replacing the ground field $K$ by $L$. Then, since the
given $K$-rational point $P$ is also defined over $L$,  we apply
Theorem \ref{thm:totally} to complete the proof.
\end{proof}
%\section*{Appendix}
%\input{try1}

\end{document}